\documentclass[12pt]{article}
\usepackage{t1enc}
\usepackage[latin1]{inputenc}
\usepackage[english]{babel}
\usepackage{amsmath,amsthm}
\usepackage{amsfonts}
\usepackage{amssymb}
\usepackage{listings}
\usepackage{enumitem}
\usepackage{lineno}
\usepackage{hyperref,cleveref}
\usepackage{latexsym}
\usepackage{graphicx}
\usepackage{epstopdf}
\usepackage{float}
\usepackage{color}
\usepackage[natural]{xcolor}
\usepackage{float,color,fancybox,shapepar,setspace,hyperref}
\theoremstyle{plain}

\newtheorem{theorem}{Theorem}

\newtheorem{claim}{Claim}
\newtheorem{lemma}[claim]{Lemma}
\numberwithin{claim}{section}

\newtheorem*{conjecture*}{Conjecture}
\newtheorem{definition}[claim]{Definition}
\newtheorem{fact}[claim]{Fact}

\newtheorem{remark}[claim]{Remark}

\newcommand{\send}[2]{c(#1 \rightarrow #2)}
\usepackage{appendix}

\title{\LARGE $4$-choosability of planar graphs with $4$-cycles far apart via the Combinatorial Nullstellensatz\thanks{This work is supported by NSFC(11971270, 11631014) of China and Shandong Province Natural Science Foundation (ZR2018MA001, ZR2019MA047) of China.}}
\author{
        Fan Yang\thanks{Data Science Institute, Shandong University, Jinan 250100, China, Email: yangfan5262@163.com.}\and
        Yue Wang\thanks{School of Mathematics and Statistics, Shandong Normal University, Jinan 250358, China, Email: wangyue\_math@163.com.}\and
        Jian-liang Wu\thanks{Corresponding author. School of Mathematics, Shandong University, Jinan 250100, China, Email: jlwu@sdu.edu.cn.}
}
\date{}

\begin{document}
\maketitle

\baselineskip 0.65cm
\begin{abstract}
By a well-known theorem of Thomassen and a planar graph depicted by Voigt, we know that every planar graph is $5$-choosable, and the bound is tight. In 1999, Lam, Xu and Liu reduced $5$ to $4$ on $C_4$-free planar graphs. In the paper, by applying the famous Combinatorial Nullstellensatz, we design an effective algorithm to deal with list coloring problems. At the same time, we prove that a planar graph $G$ is $4$-choosable if any two $4$-cycles having distance at least $5$ in $G$, which extends the result of Lam et al.
\end{abstract}


\section{Introduction}\label{section1}

All graphs considered in the paper are simple and finite. The concepts of list coloring and choosability were introduced by Vizing \cite{Viz} and independently by Erd\H{o}s, Rubin and Taylor \cite{Erd}. Given a graph $G$, a \emph{list assignment} $L$ for $G$ is a function that to each vertex $v\in V(G)$ assigns a set $L(v)$ of colors, and an \emph{$L$-coloring} is a proper coloring $\phi$ such that $\phi(v)\in L(v)$ for all $v\in V(G)$. We say that $G$ is \emph{$L$-colorable} if $G$ has an $L$-coloring. Moreover, $G$ is \emph{$k$-choosable} if $G$ is $L$-colorable for every list assignment $L$ with $|L(v)|\geq k$ for each $v \in V(G)$. List coloring is a fundamental object in graph theory with a wealth of results studying various aspects and variants. A variety of mathematicians have suggested imposing slightly stronger conditions in order to further explore the choosability of graphs, see \cite{DP, listre, 4-1-cho}. The distance of two vertices is the shortest length (number of edges) of paths between them, and the distance $d(H_1,H_2)$ of two subgraphs $H_1$ and $H_2$ is the minimum of the distances between vertices $v_1\in V(H_1)$ and $v_2\in V(H_2)$.

The classic Four Color Theorem claims that every planar graph is $4$-colorable, which was proved by Appel and Haken in 1976 \cite{4-color1, 4-color2}. However, the result can not be extended to that of list colorings as Voigt \cite{not4-ch} found a planar graph which is not $4$-choosable. Fortunately, Thomassen \cite{5-ch} proved that every planar graph is $5$-choosable by induction on the number of vertices. In order to further explore list coloring problems, forbidding certain structures within a planar graph is a common restriction used in graph coloring. Notice that all $2$-choosable graphs have been characterised by Erd\H{o}s, Rubin and Taylor \cite{Erd}. So it remains to determine whether a given planar graph is $3$- or $4$-choosable. In recent years, a number of interesting results about the choosability of special planar graphs have been obtained. Alon and Tarsi \cite{Alon-Tar} proved that every planar bipartite graph is $3$-choosable. Thomassen \cite{3-gir} showed every planar graph of girth at least $5$ is $3$-choosable, and there exist triangle-free planar graphs which are not $3$-choosable \cite{gir4not3}, so the bound $5$ is tight. Very recently, Dvo\v{r}\'{a}k \cite{Dov} showed that every planar graph in which any two $(\leq4)$-cycles have distance at least $26$ is $3$-choosable.

Steinberg's Conjecture from 1976 states that every $\{C_4, C_5\}$-free planar graph is $3$-colorable, which was disproved by Cohen-Addad et al. \cite{SteConj}. Previously, Voigt \cite{VM} disproved a list version of Steinberg's Conjecture by giving a $\{C_4, C_5\}$-free planar graph which is not $3$-choosable. A graph $G$ is said to be \emph{$k$-degenerate} if every nonempty subgraph $H$ of $G$ has a vertex of degree at most $k$ in $H$. Note that the list chromatic number of a $k$-degenerate graph is at most $k+1$. It is simple to check that every triangle-free planar graph is $3$-degenerate, and so it is $4$-choosable. In addition, it was proved that every $C_k$-free planar graph is $4$-choosable for $k=4$ in \cite{Xu}, for $k=5$ in \cite{4m, Wang-no5}, for $k=6$ in \cite{Mohar, 4m, Wang-no6}, and for $k=7$ in \cite{no7-cycle}. On the other hand, it is shown in \cite{4m} that every planar graph in which any two triangles have distance at least $2$ is $4$-choosable, and a conjecture was proposed in this paper, which claims that every planar graph without adjacent triangles is 4-choosable (this conjecture is still open so far). After that, Wang and Li \cite{Wang-no-inter3} improved one of the results in \cite{4m} by showing that each planar graph without intersecting triangles is $4$-choosable.

Inspired by the improvements of the results about triangle-free planar graphs, we further explore the picture when any two 4-cycles in a planar graph is far apart. A natural question can be proposed as follows.

\noindent
\textbf{Problem A.} Does there exist a constant $d$ such that a planar graph $G$ is $4$-choosable if any two $4$-cycles have distance at least $d$ in $G$?

We give a positive answer to this question with $d=5$.

\begin{theorem}\label{thm}
If $G$ is a planar graph such that any two $4$-cycles have distance at least $5$, then $G$ is $4$-choosable.
\end{theorem}


\section{A Structural Lemma}
For any positive integer $r$, we write $[r]$ for the set $\{1, \ldots, r\}$. Given a plane graph $G$, we denote its vertex set, edge set, face set by $V(G)$, $E(G)$, and $F(G)$, respectively. For any vertex $v\in V(G)$ (or any face $f\in F(G)$), the degree of $v$ (or $f$), denoted by $d(v)$ (or $d(f)$), is the number of edges incident with $v$ (or the length of the boundary walk of $f$, where each cut edge is counted twice). A vertex $v$ is called a $k$-vertex ($k^{+}$-vertex, or $k^{-}$-vertex) if $d(v)=k$ ($d(v)\geq k$, or $d(v)\leq k$, respectively). Analogously, a $k$-face ($k^+$-face, or $k^{-}$-face) is a face of degree $k$ (at least $k$, or at most $k$, respectively). Moreover, we use $\Delta(G)$ and $\delta(G)$ to denote the maximum degree and the minimum degree of $G$, respectively.

\begin{figure}[H]
 \begin{center}
   \includegraphics[scale=0.7]{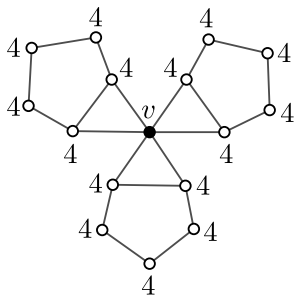}\\
   \caption[14cm]{$d(v)=6$ and $\zeta_v(f_{3b})=3$.}
   \label{zetafig}
 \end{center}
\end{figure}

We write $f=(u_1,\ldots, u_t)$ if $u_1,\ldots, u_t$ are the boundary vertices of $f$ in the clockwise order. Sometimes we replace $u_i$ with $d(u_i)$ for some $i\in[t]$ in $f=(u_1,\ldots, u_t)$ to describe the face $f$. For example, $f=(4,4,5,u_4)$ denotes a $4$-face with $d(u_1)=d(u_2)=4$, $d(u_3)=5$. For a vertex $v$ and a face $f$, let $f_k(v)$, $n_k(v)$ and $n_k(f)$ denote the number of $k$-faces incident with $v$, the number of $k$-vertices adjacent to $v$, and the number of $k$-vertices incident with $f$, respectively. Let $f=(v_1,v_2,v_3,v_4,v_5)$ be a $5$-face, $f$ is called \emph{bad} if $d(v_i)=4$ for all $i\in[5]$. For convenience, we use $f_{5b}(v)$ to denote the number of bad $5$-faces incident with a vertex $v$. In addition, let $\zeta_{v}(f_{3b})$ denote the number of $3$-faces $f=(x,y,v)$ incident with $v$ such that $d(x)=d(y)=4$ and $xy$ locates on a bad $5$-face. Below Figure \ref{zetafig} shows a $6$-vertex $v$ with $\zeta_v(f_{3b})=3$.

\noindent
A $4$-vertex $v$ with $f_3(v)+f_{5b}(v)\leq1$ of $G$ is called \emph{good}, whereas $v$ is called \emph{bad} if $f_3(v)=1$ and $f_{5b}(v)=1$.



\begin{lemma}\label{struclemma}
 Let $G$ be a connected planar graph such that any two $4$-cycles have distance at least $5$.  Then
\begin{enumerate}
  \item[$(a)$] $G$ has a $3^-$-vertex, or
  \item[$(b)$] $G$ contains one of the configurations $S_1$-$S_{47}$, see Appendix B.
\end{enumerate}
\end{lemma}

\begin{proof} Let $G$ be a counterexample to the lemma with $|V(G)|+|E(G)|$ as small as possible. Then $\delta(G)\geq 4$ and $G$ contains none of the configurations $S_1$-$S_{47}$ in Appendix B. Euler's formula $|V(G)|-|E(G)|+|F(G)|=2$ can be expressed in the form
\begin{equation}\label{veffff}
 \sum\limits_{v\in V(G)}(d_G(v)-2)+\sum\limits_{f\in F(G)}(-2)= -4.
\end{equation}
An initial charge $ch_0$ on $V(G)\cup F(G)$ is defined by letting $ch_0(v)=d(v)-2$ for each $v\in V(G)$ and $ch_0(f)=-2$ for each $f\in F(G)$. Thus we have $\sum_{z\in V(G)\cup F(G)}ch_0(z)<0$.

In the following, $\send{x}{y}$ is used to denote the amount of charges transferred from an element $x$ to an element $y$. For brevity, let $\gamma=\frac{2-\frac{1}{3}n_4(f)}{n_{5^{+}}(f)}$.

We define the following two rounds of discharging rules. The first round contains R1-R5. Let $v$ be a $k$-vertex, and let $f$ be an $\ell$-face incident with $v$.
\begin{enumerate}
 \item[$\bf R1.$] $\send{v}{f}=\frac{2}{3}$ if $\ell=3$, and $\send{v}{f}=\frac{1}{3}$ if $\ell\geq6$.

 \item[$\bf R2.$] For $k=4$ and $\ell\in\{4,5\}$.
  \begin{description}
     \item[$\bf R2.1.$] Let $T_f=\{v_i:d(v_i)=4\ \text{and}\ f_3(v_i)\le1\}$. If $f=(v_1,v_2,v_3,v_4,v_5)$ is a bad $5$-face with $f_3(v)\leq1$, then  $\send{v}{f}=\frac23$ when $|T_f|=1$, and $\send{v}{f}=\frac12$ when $|T_f|\ge2$.

     \item[$\bf R2.2.$] $\send{v}{f}=\frac{1}{3}$ otherwise.
      \end{description}

   \item[$\bf R3.$] For $k=5$, $\send{v}{f}=\frac{5}{9}$ if $\ell=4$ and $n_{6^+}(f)=1$, $\send{v}{f}=\frac{4}{9}$ if $\ell=5$ and $n_{6^+}(f)=1$, and $\send{v}{f}=\gamma$ otherwise.

 \item[$\bf R4.$] For $k\geq6$, $\send{v}{f}=\frac{7}{9}$ if $\ell=4$ and $n_{5}(f)=1$, $\send{v}{f}=\frac{5}{9}$ if $\ell=5$ and $n_{5}(f)=1$, and
       $\send{v}{f}=\gamma$ otherwise.

 \item[$\bf R5.$] Let $f=(v_1,v_2,v_3,v_4,v_5)$ be a bad $5$-face with $f_3(v_i)=2$ for each $i\in[5]$, and let $f_i=(v_i,v_{i+1},u_i)$. Then $\send{u_i}{f}=\frac{1}{9}$ if $u_i$ is not incident with any $4$-cycle.

 \end{enumerate}

Let $ch_1(x)$ be the new charge of $x$ after applying R1-R5. A vertex $v$ is called \emph{rich} if $ch_1(v)>0$ while it is called \emph{poor} if $ch_1(v)<0$ and $v$ is incident with a $4$-cycle. 
Given a poor vertex, we aim to get additional charge from rich vertices to keep it non-negative.
\begin{definition}\label{defnice}
Let $u$ be a poor vertex with  $5\leq d(u)\leq6$, and $v$ be a rich vertex. A nice $uv$-path is a path connecting $u$ and $v$ of length at most two and the internal vertex (if any) has degree at most $5$ in $G$, see Figure $\ref{npath1}$.
\end{definition}
\begin{figure}[htbp]
 \begin{center}
   \includegraphics[scale=0.7]{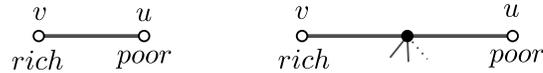}\\
   \caption[14cm]{Nice paths.}
   \label{npath1}
 \end{center}
\end{figure}
The second round R6 can be expressed as follows.

\begin{enumerate}
\item[$\bf R6.$] Let $u$ be a poor vertex, and $v_1,\ldots,v_{\ell}$ be the rich vertices at distance at most $2$ from $u$. Then $\send{v_i}{u}=ch_1(v_i)$ if $G$ has a nice $uv_i$-path.
\end{enumerate}
\begin{remark}\label{rem1}
Since the poor vertex is incident with a $4$-cycle and any two $4$-cycles have distance at least $5$, each rich vertex sends additional charge to at most one poor vertex. Note that the new charge of every rich vertex still keeps non-negative after applying R6.
\end{remark}

Let $ch_2(x)$ be the final charge of $x$ after applying R1-R6. For convenience, we say that $S_i\not\subseteq G$ if $G$ contains no subgraphs isomorphic to the configurations $S_i$ $(1\leq i\leq47)$ in Appendix B. Our goal is to show that $ch_2(z)\geq0$ for each $z\in V(G)\cup F(G)$ and so we find a contradiction to (1), which implies that the minimum counterexample does not exist. Note that $ch_2(x)=ch_1(x)$ if R6 is not applied to $x$. Thus, we have that $ch_2(f)=ch_1(f)$ for any $f\in F(G)$ by R6 and $ch_2(v)=ch_1(v)$ for any $v$ with $ch_1(v)=0$. By Remark \ref{rem1}, we get that $ch_2(v)\ge0$ for each rich vertex. So if $ch_1(z)\ge0$, then we have that  $ch_2(z)\ge0$ for each $z\in V(G)\cup F(G)$.

Since $G$ has no intersecting $4$-cycles, we immediately have the following simple fact.
\begin{fact}\label{fact1}
For each vertex $v\in V(G)$, $f_3(v)\leq\lceil\frac{d(v)}{2}\rceil$.
\end{fact}

\begin{claim}\label{c1}
 For each face $f\in F(G)$, $ch_2(f)=ch_1(f)\geq0$.
\end{claim}
\begin{proof}
If $d(f)=3$, then $ch_1(f)\geq -2+3\times \frac{2}{3}=0$ by R1. If $d(f)\geq6$, then $ch_1(f)\geq -2+6\times \frac{1}{3}=0$ by R1.

Suppose that $4\leq d(f)\leq5$ and $f$ is not a bad $5$-face. By R2.2, $f$ gets $\frac13$ from each of its incident 4-vertices.\\
(\romannumeral1) If $d(f)=4$, $n_{5}(f)=1$ and $n_{6^+}(f)=1$, then $f$ gets $\frac59$ from its incident $5$-vertex and $\frac79$ from its incident $6^+$-vertex by R3 and R4.\\
(\romannumeral2) If $d(f)=5$, $n_{5}(f)=1$ and $n_{6^+}(f)=1$, then $f$ gets $\frac49$ from its incident $5$-vertex and $\frac59$ from its incident $6^+$-vertex by R3 and R4.\\
(\romannumeral3) Otherwise, $f$ gets $\gamma$ from each of its incident $5^+$-vertices by R3 and R4.\\
Thus, we have that $ch_1(f)\geq-2+\min\{\frac13\times2+\frac59+\frac79,\frac13\times3+\frac49+\frac59,\frac13\cdot n_4(f)+\frac{2-\frac{1}{3}n_4(f)}{n_{5^{+}}(f)}\cdot n_{5^{+}}(f)\}=0$.

Suppose that $f$ is a bad $5$-face.
If there exists exactly one $i$ ($i\in[5]$) such that $f_3(v_i)\leq1$, then $f$ gets at least $\frac43$ from other incident vertices by R2.1, and so we have that $ch_1(f)\geq -2+\frac{2}{3}+4\times\frac{1}{3}=0$. If there exist at least two vertices, say $v_i$ and $v_j$, such that $f_3(v_i)\leq1$ and $f_3(v_j)\leq1$, then $f$ gets $\frac12$ from each of $v_i$ and $v_j$ by R2.1 and gets at least $1$ from other incident vertices by R2, and so we have that $ch_1(f)\geq -2+\frac12\times2+3\times\frac{1}{3}=0$. Hence, we assume that each $v_i$ satisfies $f_3(v_i)=2$. For brevity, denote by $f_i=(v_i,v_{i+1},u_i)$ the $3$-face sharing the edge $v_iv_{i+1}$ with $f$, and let $U=\{u_1,u_2,u_3,u_4,u_5\}$. Since $S_2\nsubseteq G$, we get that $d(u_i)\geq5$ for each $i\in[5]$. By the assumption of $G$, either at most one vertex in $U$ lies on a $4$-cycle, or two vertices in $U$ lie on the same $4$-cycle. Let $U^{*}\subseteq U$ such that each vertex in $U^{*}$ does not lie on any $4$-cycle. Note that $|U^{*}|\geq3$ and it follows that $f$ receives at least $3\times\frac19$ from $U^{\ast}$ by R5. So we get that $ch_2(f)=ch_1(f)\geq-2+5\times\frac{1}{3}+3\times\frac{1}{9}=0$ by R2.2.
\end{proof}


\begin{claim}\label{c2}
 For each $4$-vertex $v$, $ch_1(v)\geq 0$.
 In particular, for each good $4$-vertex $v$, $ch_1(v)\ge\frac13$.
\end{claim}
\begin{proof}
Let $v$ be a $4$-vertex. If $f_3(v)+f_{5b}(v)\leq2$, then $ch_1(v)\ge2-\frac23\times2-\frac13\times2=0$ by R1 and R2. So suppose that $f_3(v)+f_{5b}(v)\ge3$. By Fact \ref{fact1}, $f_3(v)\le2$. As $S_2,S_3\not\subseteq G$, we have that $f_{5b}(v)\le2$ and if $f_{5b}(v)=2$, then $f_3(v)=0$. It remains to consider the case that $f_3(v)=2$ and $f_{5b}(v)$=1. By R2, $v$ sends $\frac13$ to each of other $4^+$-faces and $\frac23$ to each 3-face. Thus, $ch_1(v)\ge2-\frac23\times2-\frac13\times2=0$.

Since $f_3(v)+f_{5b}(v)\le1$ holds for each good vertex $v$,  we have that $ch_1(v)\ge2-\frac13\times3-\frac23=\frac13$ by R1-R2.
\end{proof}
\begin{claim}\label{c3}
$ch_1(v)\geq 0$ if $v$ is a $7^+$-vertex, or a $6$-vertex with $f_{6^+}(v)\geq1$, or a $5$-vertex with $f_{6^+}(v)\geq2$.
\end{claim}
\begin{proof}
Let $v$ be a vertex. Suppose $d(v)$ is odd. Note that $f_3(v)\leq\frac{d(v)+1}{2}$ by Fact \ref{fact1}. If $f_3(v)=\frac{d(v)+1}{2}$, then by R1 and R4, we have that $ch_1(v)\geq d(v)-2-\frac{2}{3}d(v)=\frac{d(v)-6}{3}$.
If $f_3(v)\leq\frac{d(v)-1}{2}$ and $f_4(v)=1$, then by R1 and R4, we have that $ch_1(v)\geq d(v)-2-1-\frac{2}{3}\left(d(v)-1\right) =\frac{d(v)-7}{3}$. If $f_3(v)\leq\frac{d(v)-1}{2}$ and $f_4(v)=0$, then by R1 and R4-R5, we have that $ch_1(v)\geq d(v)-2-\frac{2}{3}d(v)-\frac{1}{9}\left(\frac{d(v)-1}{2}\right) =\frac{5d(v)-35}{18}$.
Particularly, if $f_{6^+}(v)\geq2$, then $ch_1(v)\geq\min\left\{\frac{d(v)-7}{3}, \frac{5d(v)-35}{18}\right\}+2\times\frac{1}{3}=\min\left\{\frac{d(v)-5}{3}, \frac{5d(v)-23}{18}\right\}$.



Suppose $d(v)$ is even. Note that $f_3(v)\leq \frac{d(v)}{2}$ by Fact \ref{fact1}. If $f_4(v)=1$, then by R1 and R4, we have that $ch_1(v)\geq d(v)-2-1-\frac{2}{3}\left(d(v)-1\right)  =\frac{d(v)-7}{3}.$ If $f_4(v)=0$, then by R1 and R4-R5, we have that  $ch_1(v)\geq d(v)-2-\frac{2}{3}d(v)-\frac{1}{9}\left(\frac{d(v)}{2}\right)=\frac{5d(v)-36}{18}.$ In particular, if $f_{6^+}(v)\geq1$, then $ch_1(v)\geq\min\left\{\frac{d(v)-7}{3},\frac{5d(v)-36}{18}\right\}+\frac13=\min\left\{\frac{d(v)-6}{3},\frac{5d(v)-30}{18}\right\}$.

Therefore, Claim \ref{c3} is true.
\end{proof}
Now it remains to consider the vertices of $W_1=\{v: d(v)=6 \ and\  f_{6^+}(v)=0\}$ and $W_2=\{v: d(v)=5 \ and\  f_{6^+}(v)\le1\}$ by Claim \ref{c2} and \ref{c3}. 

For $v\in W_1$, let $N(v)=\{v_1,\ldots,v_6\}$ and let $f_1, \ldots, f_6$ be the faces incident with $v$ in clockwise such that $v_i$ and $v_{i+1}$ are incident with $f_i$. In the following Claims \ref{W11}-\ref{W14}, we show that $ch_2(v)\geq0$ for each vertex $v\in W_1$.

\begin{claim}\label{W11}
For each vertex $v\in W_1$ with $f_3(v)\le2$ and  $f_4(v)=1$, $ch_2(v)\geq0$.
\end{claim}
\begin{proof}
W.l.o.g., let $f_1$ be the $4$-face, denoted by $v_1vv_2x$. Note that $v$ sends no charge to a bad $5$-face (if it exists) which is incident with a $(4,4,v)$-face by R5. According to R1 and R4, $v$ sends at most $1$ to each 4-face and $\frac23$ to each 3-face and 5-face. Thus, $ch_1(v)\ge4-1-\frac23\times5=-\frac13$. If $ch_1(v)\ge0$, then we are done. So suppose that $ch_1(v)<0$, that is, $v$ is poor. Clearly, if there is a good 4-vertex in $N(v)$, then  $ch_2(v)\geq4-1-\frac23\times5+\frac13=0$ by Claim \ref{c2} and R6. Next we only consider the case that there is no good 4-vertex in $N(v)$.

Now we first claim that $f_i$ is not a $(4,4,4,4,6)$-face for each $i\in\{2,6\}$ (that is, $n_{5^+}(f_i)\geq2$). Suppose to the contrary that for some $i\in\{2,6\}$, $f_i$ is a $(4,4,4,4,6)$-face, say $f_2$. As $S_2, S_3\not\subseteq G$, we get that $f_3(v_2)+f_{5b}(v_2)\le1$ and $v_2$ is a good 4-vertex, a contradiction. Similarly, if $f_6$ is a $(4,4,4,4,6)$-face, then $v_1$ is a good 4-vertex, a contradiction.
\begin{figure}[H]
 \begin{center}
  \includegraphics[scale=0.7]{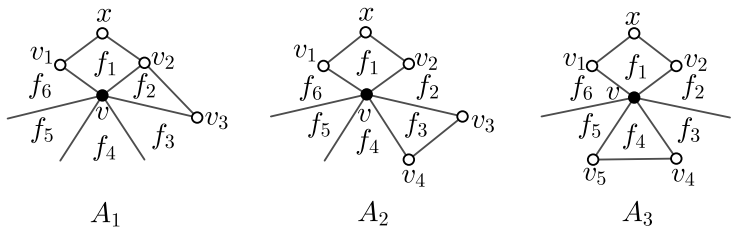}
  \caption[14cm]{Configurations for $6$-vertex $v$ with $f_4(v)=1$.}
   \label{6-vx2eps}
 \end{center}
\end{figure}
 \textbf{Case 1.}  $n_{5^+}(f_1)=1$.

 \textbf{Subcase 1.1.} Assume that $f_3(v)\le1$. We will show that there are at least three $5$-faces $f_i$ such that $n_{5^+}(f_i)\geq2$, which implies that $\send{v}{f_i}\leq \frac{5}{9}$ by R4, and so $ch_1(v)\geq4-1-2\times\frac{2}{3}-3\times\frac{5}{9}=0$ by R1.

 $(a)$ Suppose that $f_3(v)=0$. Since $G$ has no intersecting $4$-cycles, the remaining faces incident with $v$ are all $5$-faces. By $S_{24}\not\subseteq G$, there exists at least one $i$ ($i\in\{3,4,5\}$) such that $n_{5^+}(f_i)\geq2$. Note that $n_{5^+}(f_2)\geq2$ and $n_{5^+}(f_6)\geq2$, so we are done.

 $(b)$ Suppose that $f_3(v)=1$. By symmetry, three cases need to be considered (see Figure \ref{6-vx2eps}). In $A_1$, since $S_{32}\not\subseteq G$, we have that $n_{5^+}(f_3)\geq2$. In $A_2$, since $S_{24}\not\subseteq G$, we have that $n_{5^+}(f_4)\geq2$ or $n_{5^+}(f_5)\geq2$. In $A_3$, since  $S_{27}\not\subseteq G$, we have that $n_{5^+}(f_3)\geq2$ or $n_{5^+}(f_5)\geq2$. Note that if $f_i$ is a 5-face, then $n_{5^+}(f_i)\geq2$ for $i\in\{2,6\}$, so we are done.

 \textbf{Subcase 1.2.} Assume that $f_3(v)=2$. There are four subcases to be considered.

Firstly, we suppose that $d(f_2)=d(f_4)=3$. Note that $f_3(v_1)\leq1$ and $v_1$ is not good. It implies that $v_1$ is bad. Since $S_2,S_3\not\subseteq G$, $v_1x$ locates on the same bad $5$-face. In this situation, $f_3(x)\leq1$, and by R2.1, each of $\{v_1, x\}$ sends $\frac12$ to the bad 5-face. Thus, by R2 $ch_1(u)\ge2-\frac23-\frac12-\frac13\times2=\frac16$ for each $u\in\{v_1,x\}$. Therefore, each of $\{v_1, x\}$ sends $\frac16$ to $v$ (if $ch_1(v)<0$) via a nice path by R6. Thus $ch_2(v)\geq0$. The case that $d(f_2)=d(f_5)=3$ is similar as above.


Next, we suppose that $d(f_3)=d(f_5)=3$. Since $v_1$ and $v_2$ are not good and $S_2,S_3\not\subseteq G$, $v_1x$ locates on the same bad $5$-face $g_1$ and $v_2x$ locates on the same bad $5$-face $g_2$. By $S_2\not\subseteq G$, we have that $f_3(x)=0$. Note that $f_3(v_i)\le1$ for each $i\in[2]$. It follows that $|T{g_1}|\ge2$ and $|T{g_2}|\ge2$.   Thus,
$ch_1(x)\geq 2-2\times \frac{1}{2}-2\times \frac{1}{3}=\frac{1}{3}$ by R1 and R2. Hence, $v$ (if $ch_1(v)<0$) could receive at least $\frac{1}{3}$ from $x$ via a nice path by R6, and $ch_2(v)\geq0$.

It remains to consider the case where $d(f_2)=d(f_6)=3$. Since $S_{32}\nsubseteq G$, we get that for each $i\in\{3,5\}$, $n_{5^+}(f_i)\geq2$ and $\send{v}{f_i}\leq \frac{5}{9}$ by R4. If $n_{5^+}(f_4)\geq2$, then $ch_1(v)\geq 4-1-2\times \frac{2}{3}-3\times \frac{5}{9}=0$ by R1 and R4. Now let $n_{5^+}(f_4)=1$, and denote by $f_4=(v,v_4,y_1,y_2,v_5)$, that is $d(v_4)=d(v_5)=d(y_1)=d(y_2)=4$. Note that $f_3(v_4)=f_3(v_5)\leq1$. So we may assume that both $v_4$ and $v_5$ are not good (otherwise $v$ receives at least $\frac{1}{3}$ from $\{v_4,v_5\}$ and $ch_2(v)\geq0$). Since $S_2\nsubseteq G$, $v_4y_1$ and $v_5y_2$ locate on two bad $5$-faces, respectively. On the other hand, notice that $S_2, S_{47}\nsubseteq G$, and then at least one $j\in\{3,5\}$ satisfying $n_{5^+}(f_j)\geq 3$, and so $\send{v}{f_j}\leq \frac{4}{9}$ for some $j\in\{3,5\}$ by R4. Thus $ch_1(v)\geq4-1-3\times \frac{2}{3}-\frac{5}{9}-\frac{4}{9}=0$ by R1 and R4.

 \textbf{Case 2.} $n_{5^+}(f_1)\geq2$. Since $S_2, S_{36} \not\subseteq G$, there exists at least one $i\in\{2,3,4,5,6\}$ such that $n_{5^+}(f_i)\geq2$, and we have $\send{v}{f_{i}}\leq \frac{5}{9}$ by R4. So $ch_1(v)\geq4-4\times\frac{2}{3}-\frac{7}{9}-\frac{5}{9}=0$ by R1 and R4.
\end{proof}
\begin{claim}\label{W12}
For each vertex $v\in W_1$ with $f_3(v)\le2$ and  $f_4(v)=0$, $ch_2(v)\geq0$.
\end{claim}
\begin{proof}
Suppose that $f_4(v)=0$. If $f_3(v)=0$, then $ch_2(v)\geq4-6\times\frac{2}{3}=0$ by R1 and R4. If $f_3(v)=1$, then by $S_{35}\nsubseteq G$, either $\zeta_v(f_{3b})=0$ or there exists some $i$ such that $n_{5^+}(f_i)\geq2$, and thus $ch_1(v)\geq4-5\times\frac{2}{3}-\max\{\frac{5}{9}+\frac{1}{9}, \frac{2}{3}\}=0$ by R1, R4-R5. Finally, we discuss the case where $f_3(v)=2$. If the two $3$-faces are consecutive, then $ch_1(v)\geq4-6\times\frac{2}{3}=0$ by R1 and R4. Otherwise if they are not consecutive, by the fact that $S_{30},S_{35}\nsubseteq G$, we get that $ch_1(v)\geq4-2\times\frac{5}{9}-4\times\frac{2}{3}-2\times\frac{1}{9}=0$ by R1, R4-R5.
\end{proof}

Next we focus on the case $f_3(v)=3$. Since $S_{46}\nsubseteq G$, we get $\zeta_v(f_{3b})\leq2$.
\begin{figure}[H]
 \begin{center}
  \includegraphics[scale=0.7]{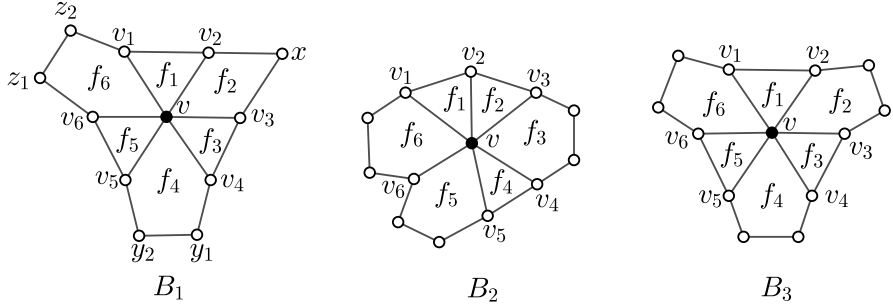}
  \caption[14cm]{Configurations for $6$-vertex $v$ with $f_3(v)=3$.}
   \label{6-4feps}
 \end{center}
\end{figure}
Recall that $v$ sends no charge to a bad $5$-face which is incident with a $(4,4,v)$-face by R5.
\begin{claim}\label{W13}
For each vertex $v\in W_1$ with $f_3(v)=3$ and $f_4(v)=1$ $($see Figure $\ref{6-4feps}$$(B_1)$$)$, $ch_2(v)\geq0$.
\end{claim}
\begin{proof}
We divide the proof into four possibilities depending on $n_4(v)\in\{4,5,6\}$ or $n_4(v)\leq3$.

(\romannumeral1) Suppose $n_4(v)=6$.
$(a)$. $d(x)=4$. As $S_{32}\nsubseteq G$, for each $i\in\{4,6\}$, we have $n_{5^+}(f_i)\geq2$, and so $\send{v}{f_{i}}\leq \frac{5}{9}$ by R4. If $n_{6^+}(f_i)\geq2$ or $n_{5^+}(f_i)\geq3$ for each $i\in\{4,6\}$, then $ch_1(v)\geq4-3\times\frac{2}{3}-1-2\times\frac{1}{2}=0$ by R1 and R4. Assume $n_{6}(f_4)=1$ and $n_{5^+}(f_4)=2$. Denote by $f_4=(v,v_4,y_1,y_2,v_5)$.
First, let $d(y_1)=5$ and $d(y_2)=4$. If $f_3(v_4)=1$, then according to $S_{25}\nsubseteq G$, $v_4$ can not locate on a bad $5$-face. Thus $v_4$ is good, and $ch_1(v_4)\geq\frac{1}{3}$. If $f_3(v_4)=2$,
then by $S_{25},S_{41},S_{44}\nsubseteq G$, we have $\zeta_{y_1}(f_{3b})=0$. Thus $ch_1(y_1)\geq3-3\times\frac{2}{3}-\frac{4}{9}-\frac{1}{2}=\frac{1}{18}$ by R1 and R3. In both cases, $\{v_4,y_1\}$ could send at least $\frac{1}{18}$ to $v$ (if $ch_1(v)<0$) via a nice path by R6. Second, let $d(y_1)=4$ and $d(y_2)=5$. If $f_3(v_5)=1$, then $v_5$ can not locate on a bad $5$-face by $S_{38}\nsubseteq G$. Thus $v_5$ is good, and $ch_1(v_5)\geq\frac{1}{3}$. If $f_3(v_5)=2$,
then by $S_{38},S_{39},S_{42}\nsubseteq G$, we have $\zeta_{y_2}(f_{3b})=0$. Thus $ch_1(y_2)\geq3-3\times\frac{2}{3}-\frac{4}{9}-\frac{1}{2}=\frac{1}{18}$ by R1 and R3. In both cases, $\{v_5,y_2\}$ could send at least $\frac{1}{18}$ to $v$ (if $ch_1(v)<0$) via a nice path by R6. In conclusion, $v$ could receive at least $\frac{1}{18}$ from $\{v_4,v_5,y_1,y_2\}$. By symmetry, the same arguments also hold for the vertices on $f_6$ (i.e. $\{v_1,v_6,z_1,z_2\}$). If $n_{6}(f_6)=1$ and $n_{5^+}(f_6)=2$, then $ch_2(v)\geq4-3\times\frac{2}{3}-1-2\times\frac{5}{9}+2\times\frac{1}{18}=0$ by R1, R4 and R6.
Otherwise $ch_2(v)\geq4-3\times\frac{2}{3}-1-\frac{5}{9}-\frac{1}{2}+\frac{1}{18}=0$.

$(b)$. $d(x)\geq5$.
Since $S_{30}\nsubseteq G$, $n_{5^+}(f_i)\geq2$ for some $i\in\{4,6\}$, and $\send{v}{f_{i}}\leq \frac{5}{9}$ by R4. Thus $ch_1(v)\geq4-4\times\frac{2}{3}-\frac{7}{9}-\frac{5}{9}=0$ by R1 and R4.

(\romannumeral2) Suppose $n_4(v)=5$. By symmetry, we only need to consider three subcases: $d(v_1)\geq5$, $d(v_2)\geq5$ and $d(v_5)\geq5$.

$(a)$. $d(x)=4$.
Assume that $d(v_1)\geq5$. If $d(v_1)=5$, then $n_{5^+}(f_6)\geq3$ by $S_{33}\not\subseteq G$, and we have $\send{v}{f_{i}}\leq \frac{4}{9}$ by R4. Thus $ch_1(v)\geq4-3\times\frac{2}{3}-1-\frac{5}{9}-\frac{4}{9}=0$ by R1 and R4. If $d(v_1)=6$, then $ch_1(v_1)\geq4-4\times\frac{2}{3}-\frac{1}{2}-\frac{5}{9}-2\times\frac{1}{9}=\frac{1}{18}$ by R1, R4-R5 because of $S_{40}\nsubseteq G$. If $d(v_1)\geq7$, then by Claim \ref{c3}, $ch_1(v_1)\geq \frac{5d(v)-36}{18}+
\frac{1}{9}+\frac{1}{18}\geq\frac{1}{18}$. Hence, when $d(v_1)\geq6$, $v$ (if $ch_1(v)<0$) could receive at least $\frac{1}{18}$ from $v_1$ via a nice path by R6. Thus $ch_2(v)\geq4-3\times\frac{2}{3}-1-\frac{5}{9}-\frac{1}{2}+\frac{1}{18}=0$ by R1, R4 and R6.

Assume that $d(v_2)\geq5$. Then $n_{5^+}(f_i)\geq2$ for some $i\in\{4,6\}$ by $S_{30}\nsubseteq G$, and we have $\send{v}{f_{i}}\leq \frac{5}{9}$ by R4. Thus $ch_1(v)\geq4-4\times\frac{2}{3}-\frac{7}{9}-\frac{5}{9}=0$ by R1 and R4.

Assume that $d(v_5)\geq5$. Since $S_{32}\not\subseteq G$, $n_{5^+}(f_6)\geq2$ and $\send{v}{f_{6}}\leq \frac{5}{9}$ by R4. According to $S_{26}\nsubseteq G$, either $d(v_5)\geq6$ or $n_{5^+}(f_4)\geq3$. If $n_{5^+}(f_4)\geq3$, then $\send{v}{f_{4}}\leq \frac{4}{9}$ by R4, and thus $ch_1(v)\geq4-3\times\frac{2}{3}-1-\frac{5}{9}-\frac{4}{9}=0$ by R1 and R4; if $d(v_5)\geq6$, then by the similar arguments as above, we have that $v$ (if $ch_1(v)<0$) could receive at least $\frac{1}{18}$ from $v_5$ via a nice path by R6, and thus $ch_2(v)\geq4-3\times\frac{2}{3}-1-\frac{5}{9}-\frac{1}{2}+\frac{1}{18}=0$ by R1, R4 and R6.

$(b)$. $d(x)\geq5$.
In all three cases, it is easy to check that $ch_1(v)\geq4-4\times\frac{2}{3}-\max\{\frac{7}{9}+\frac{5}{9}, \frac{5}{9}+\frac{2}{3}\}=0$ by R1 and R4.

(\romannumeral3) Suppose $n_4(v)=4$. That is, $n_{5^+}(v)=2$. If the pair of two $5^+$-vertices fall in $\{(v_1,v_2), (v_1,v_3), (v_2,v_3), (v_2,v_5), (v_2,v_6)\}$, then we have $ch_1(v)\geq4-4\times\frac{2}{3}-\max\{\frac{7}{9}+\frac{5}{9}, \frac{2}{3}+\frac{5}{9}\}=0$ by R1 and R4. By symmetry, it remains to discuss the following cases.

Assume that $d(v_1)\geq5$ and $d(v_4)\geq5$.
$(a)$. $d(x)=4$. Note that $v_1$ and $v_4$ are symmetric to some extent. If $d(v_1)\geq6$ and $d(v_4)\geq6$, then $ch_1(v)\geq4-3\times\frac{2}{3}-1-2\times\frac{1}{2}=0$ by R1 and R4. If $d(v_i)=5$ for some $i\in\{1,4\}$, then  $n_{5^+}(f_{8-2i})\geq3$ by $S_{31}\not\subseteq G$, and so $\send{v}{f_{8-2i}}\leq \frac{4}{9}$ by R4. Thus $ch_1(v)\geq4-3\times\frac{2}{3}-1-\max\{\frac{1}{2}+\frac{4}{9},2\times\frac{4}{9}\}=\frac{1}{18}>0$ by R1 and R4. $(b)$. $d(x)\geq5$. Then $ch_1(v)\geq4-3\times\frac{2}{3}-\frac{7}{9}-2\times\frac{5}{9}=\frac{1}{9}>0$ by R1 and R4.

Assume that $d(v_1)\geq5$ and $d(v_5)\geq5$.
$(a)$. $d(x)=4$. If $d(v_1)\geq6$ and $d(v_5)\geq6$, then $ch_1(v)\geq4-3\times\frac{2}{3}-1-2\times\frac{1}{2}=0$ by R1 and R4. If $d(v_i)=5$ for some $i\in\{1,5\}$, then $n_{5^+}(f_{\frac{13-i}{2}})\geq3$ by $S_{26}, S_{33}\not\subseteq G$, and so $\send{v}{f_{\frac{13-i}{2}}}\leq \frac{4}{9}$ by R4. Thus $ch_1(v)\geq4-3\times\frac{2}{3}-1-\max\{\frac{1}{2}+\frac{4}{9}, 2\times\frac{4}{9}\}=\frac{1}{18}>0$ by R1 and R4.
$(b)$. $d(x)\geq5$. Then $ch_1(v)\geq4-3\times\frac{2}{3}-2\times\frac{5}{9}-\frac{7}{9}=\frac{1}{9}>0$ by R1 and R4.

Assume that $d(v_1)\geq5$ and $d(v_6)\geq5$.
$(a)$. $d(x)=4$. Since $S_{32}\nsubseteq G$, we get $n_{5^+}(f_4)\geq2$, and $\send{v}{f_{4}}\leq \frac{5}{9}$ by R4. Then we have $ch_1(v)\geq4-3\times\frac{2}{3}-1-\frac{5}{9}-\frac{4}{9}=0$ by R1 and R4.
$(b)$. $d(x)\geq5$. Then we have $ch_1(v)\geq4-4\times\frac{2}{3}-\frac{7}{9}-\frac{4}{9}=\frac{1}{9}>0$ by R1 and R4.

Assume that $d(v_5)\geq5$ and $d(v_6)\geq5$. $(a)$. $d(x)=4$. If $d(v_5)\geq6$ and $d(v_6)\geq6$, then $ch_1(v)\geq4-3\times\frac{2}{3}-1-2\times\frac{1}{2}=0$ by R1 and R4. If $d(v_i)=5$ and $d(v_{11-i})\geq6$ for some $i\in\{5,6\}$, then by $S_2\nsubseteq G$, we get that $ch_1(v_i)\geq3-\frac{2}{3}-2\times\frac{4}{9}-\max\{\frac{2}{3}+\frac{1}{2}+\frac{1}{9}, 2\times\frac{2}{3}\}=\frac{1}{9}$ by R1, R3 and R5. Hence, $v_i$ could send at least $\frac{1}{9}$ to $v$ via a nice path by R6, and $ch_2(v)\geq4-3\times\frac{2}{3}-1-\frac{5}{9}-\frac{1}{2}+\frac{1}{9}>0$ by R1, R4 and R6. If $d(v_5)=5$ and $d(v_6)=5$, then there is at least one $i\in\{4,6\}$ such that $n_{5^+}(f_i)\geq3$ by $S_{34}\not\subseteq G$, and so $\send{v}{f_{i}}\leq \frac{4}{9}$ by R4. Hence, $ch_1(v)\geq4-3\times\frac{2}{3}-1-\frac{5}{9}-\frac{4}{9}=0$ by R1 and R4. $(b)$. $d(x)\geq5$. Then $ch_1(v)\geq4-3\times\frac{2}{3}-\frac{7}{9}-2\times\frac{5}{9}=\frac{1}{9}>0$ by R1 and R4.

(\romannumeral4) Suppose $n_4(v)\leq3$. If $n_4(v)=3$, then we have $ch_1(v)\geq4-3\times\frac{2}{3}-\max\{2\times\frac{5}{9}+\frac{7}{9}, \frac{2}{3}+\frac{4}{9}+\frac{7}{9}, \frac{2}{3}+2\times\frac{5}{9}, 1+\frac{4}{9}+\frac{5}{9}, \frac{4}{9}+\frac{5}{9}+\frac{7}{9}\}=0$ by R1 and R4. If  $n_4(v)=2$, then we have $ch_1(v)\geq4-3\times\frac{2}{3}-\max\{\frac{7}{9}+\frac{4}{9}+\frac{5}{9}, 3\times\frac{5}{9}, \frac{2}{3}+\frac{4}{9}+\frac{5}{9}, 1+2\times\frac{4}{9}\}=\frac{1}{9}>0$ by R1 and R4. If $n_4(v)=1$, then we have $ch_1(v)\geq4-3\times\frac{2}{3}-\max\{2\times\frac{5}{9}+\frac{4}{9},2\times\frac{4}{9}+\frac{7}{9}\}=\frac{4}{9}>0$ by R1 and R4. If $n_4(v)=0$, then we have $ch_1(v)\geq4-3\times\frac{2}{3}-\frac{5}{9}-2\times\frac{4}{9}=\frac{5}{9}>0$ by R1 and R4.
\end{proof}
\begin{claim}\label{W14}
For each vertex $v\in W_1$ with $f_3(v)=3$ and $f_4(v)=0$, $ch_2(v)\geq0$.
\end{claim}
\begin{proof}
If $v$ is incident with a $4$-cycle (see Figure \ref{6-4feps}$(B_2)$), then $v$ also sends no charge to a bad $5$-face (if it exists) which is incident with a $(4,4,v)$-face by R5. Thus $ch_1(v)\geq4-6\times\frac{2}{3}=0$ by R1 and R4. Next we turn to the case that $v$ is not incident with any $4$-cycle, see Figure \ref{6-4feps}$(B_3)$. Recall that $\zeta_v(f_{3b})\leq2$.

(\romannumeral1) Suppose $n_4(v)=6$. Then there are at least two faces $f_i$, $f_j$ in $\{f_2,f_4,f_6\}$ satisfying $f_i \neq(4,4,4,4,6)$ and $f_j\neq (4,4,4,4,6)$ by $S_{30}\nsubseteq G$. Thus $ch_1(v)\geq4-4\times\frac{2}{3}-2\times\frac{5}{9}-2\times\frac{1}{9}=0$ by R1 and R4-R5.

(\romannumeral2) Suppose $n_4(v)=5$. By symmetry, say $d(v_1)\geq5$. Since $S_{30}\nsubseteq G$, we get $n_{5^+}(f_i)\geq2$ when $d(f_i)=5$ for some $i\in\{2,4\}$, and so $\send{v}{f_{i}}\leq \frac{5}{9}$ by R4. Thus $ch_1(v)\geq4-4\times\frac{2}{3}-2\times\frac{5}{9}-2\times\frac{1}{9}=0$ by R1 and R4-R5.

(\romannumeral3) Suppose $n_4(v)\leq4$. If $n_4(v)=4$, then $ch_1(v)\geq4-\max\{4\times\frac{2}{3}+2\times\frac{5}{9}+2\times\frac{1}{9}, 5\times\frac{2}{3}+\frac{4}{9}+\frac{1}{9}\}=0$ by R1 and R4-R5. If $n_4(v)=3$, then  $ch_1(v)\geq4-\max\{4\times\frac{2}{3}+\frac{4}{9}+\frac{5}{9}+\frac{1}{9}, 3\times\frac{2}{3}+3\times\frac{5}{9}+\frac{1}{9}\}=\frac{2}{9}>0$ by R1, R4-R5. If $n_4(v)=2$, then $ch_1(v)\geq4-\max\{3\times\frac{2}{3}+\frac{4}{9}+2\times\frac{5}{9}+\frac{1}{9}, 4\times\frac{2}{3}+2\times\frac{5}{9}\}=\frac{2}{9}>0$ by R1, R4-R5. If $n_4(v)=1$, then we have $ch_1(v)\geq4-3\times\frac{2}{3}-2\times\frac{4}{9}-\frac{5}{9}=\frac{5}{9}>0$ by R1 and R4. If $n_4(v)=0$,  then we have $ch_1(v)\geq4-3\times\frac{2}{3}-3\times\frac{4}{9}=\frac{2}{3}>0$ by R1 and R4.
\end{proof}
For each vertex $v\in W_2$, denote by $f_i$ $(i\in[5])$ the faces incident with $v$. If $d(f_i)=3$ for some $i$, then denote by $f_i=(v,v_i,v_{i+1})$. The following Claims \ref{W21}-\ref{W25} imply that  $ch_2(v)\geq0$, for each vertex $v\in W_2$.

\begin{claim}\label{W21}
For each vertex $v\in W_2$ with $f_3(v)=3$, $ch_2(v)\ge0$.
\end{claim}
\begin{proof}
In this case, $f_4(v)=0$ since $G$ does not contain intersecting $4$-cycles. Let $f_1$, $f_2$ and $f_4$ be the $3$-faces incident with $v$. If $d(f_i)\geq6$ for some $i\in\{3,5\}$, then  $ch_1(v)\geq3-4\times\frac{2}{3}-\frac{1}{3}=0$ by R1 and R3. Next, we consider the situation where $d(f_i)=5$ for each $i\in\{3,5\}$, see Figure \ref{fv3eps}$(C_1)$.
\begin{figure}[H]
 \begin{center}
   \includegraphics[scale=0.7]{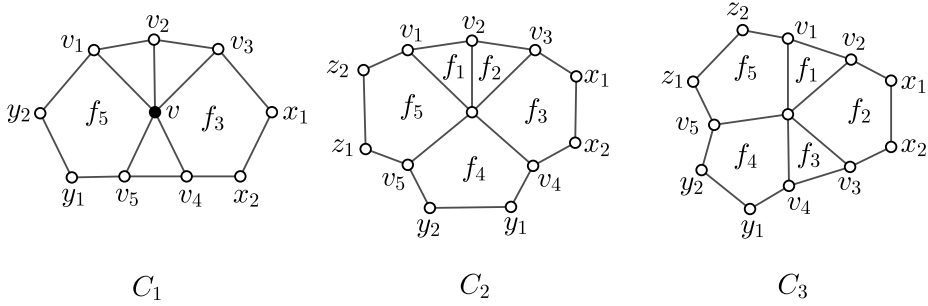}\\
   \caption[14cm]{Configurations for $5$-vertex $v$.}
   \label{fv3eps}
 \end{center}
\end{figure}
(\romannumeral1) Suppose $n_4(v)=5$. Then $n_{5^+}(f_3)\geq2$ and $n_{5^+}(f_5)\geq2$ hold by $S_2\nsubseteq G$, and so $\send{v}{f_{i}}\leq \frac{1}{2}$ for each $i\in\{3,5\}$ by R3. Thus $ch_1(v)\geq3-3\times\frac{2}{3}-2\times\frac{1}{2}=0$ by R1 and R3.

(\romannumeral2) Suppose $n_4(v)=4$, that is $n_{5^+}(v)=1$. By symmetry, there are only three cases need to be considered: $d(v_1)\geq5$; $d(v_2)\geq5$; $d(v_4)\geq5$. In all three cases, since $S_{2}\not\subseteq G$, we have $ch_1(v)\geq3-3\times\frac{2}{3}-2\times\frac{1}{2}=0$ by R1 and R3.

(\romannumeral3) Suppose $n_4(v)=3$, that is $n_{5^+}(v)=2$. If the pair of two $5^{+}$-vertices fall in $\{(v_1,v_2),(v_1,v_3), (v_1,v_4),(v_1,v_5),(v_4,v_5)\}$, then $ch_1(v)\geq3-3\times\frac{2}{3}-\max\{2\times\frac{1}{2}, \frac{4}{9}+\frac{1}{2}\}=0$ by R1 and R3. By symmetry, it remains to consider the pair $(v_2,v_4)$ with $d(v_2)\geq5$ and $d(v_4)\geq5$. We may assume that $n_{5^+}(f_5)=1$ (otherwise $ch_1(v)\geq0$).
Since $S_{22}\not\subseteq G$, $d(v_2)\geq6$. If $d(v_2)=6$, then by $S_{29},S_{35},S_{37}\nsubseteq G$, there are at least two faces $\tilde{f}$ incident with $v_2$ such that $n_{5^+}(\tilde{f})\geq2$, and so $\send{v}{\tilde{f}}\leq \frac{5}{9}$. Thus $ch_1(v_2)\geq4-4\times\frac{2}{3}-2\times\frac{5}{9}=\frac{2}{9}$ by R1 and R4. If $d(v_2)\geq7$, then $ch_1(v_2)\geq\frac{d(v)-6}{3}>\frac{2}{9}$. Hence, $v_2$ could send at least $\frac{2}{9}$ to $v$ (if $ch_1(v)<0$) via a nice path by R6, and $ch_2(v)\geq3-4\times\frac{2}{3}-\frac{1}{2}+\frac{2}{9}>0$ by R1, R3 and R6.

(\romannumeral4) Suppose $n_4(v)=2$. If the pair of two $4$-vertices fall in $\{(v_1,v_2),$ $(v_1,v_3), (v_1,v_4),(v_2,v_4),$ $(v_4,v_5)\}$, then $ch_1(v)\geq3-3\times\frac{2}{3}-\max\{2\times\frac{1}{2}, \frac{4}{9}+\frac{1}{2}\}=0$ by R1 and R3. By symmetry, it remains to consider the pair $(v_1,v_5)$ with $d(v_1)=d(v_5)=4$. We may assume that $n_{5^+}(f_5)=1$ (otherwise $ch_1(v)\geq0$). If $d(v_2)=6$, then $ch_1(v)\geq4-5\times\frac{2}{3}-\frac{5}{9}=\frac{1}{9}$ by R1 and R4; if $d(v_2)\geq7$, then $ch_1(v_2)\geq\frac{d(v)-6}{3}>\frac{1}{9}$. Hence, $v_2$ could send at least $\frac{1}{9}$ to $v$ (if $ch_1(v)<0$) via a nice path by R6, and $ch_2(v)\geq3-4\times\frac{2}{3}-\frac{4}{9}+\frac{1}{9}=0$ by R1, R3 and R6. The same results hold for $v_3$. We now turn to the case $d(v_2)=d(v_3)=5$. For simplicity, denote by $f_6$, $f_7$ and $f_8$ the remaining faces incident with $v_3$ in clockwise. If $n_{5^+}(f_6)\geq3$, then $\send{v_3}{f_{6}}\leq \frac{4}{9}$ by R4, and thus  $ch_1(v_3)\geq3-3\times\frac{2}{3}-2\times\frac{4}{9}=\frac{1}{9}$ by R1 and R4. Otherwise, $n_{5^+}(f_6)=2$. If $d(f_7)=3$, then $n_{5^+}(f_8)\geq2$ as $S_2\nsubseteq G$; if $d(f_8)=3$, then by $S_{20}$, $n_{5^+}(f_7)\geq2$; if none of $f_7$ and $f_8$ are $3$-faces, then by $S_3$, $n_{5^+}(f_i)\geq2$ for some $i\in\{7,8\}$.
In all cases, we have $ch_1(v_3)\geq3-2\times\frac{2}{3}-2\times\frac{1}{2}-\frac{4}{9}=\frac{2}{9}$ by R1 and R3. Thus $v_3$ could send at least $\frac{2}{9}$ to $v$ (if $ch_1(v)<0$) via a nice path by R6, and $ch_2(v)\geq3-4\times\frac{2}{3}-\frac{4}{9}+\frac{2}{9}>0$ by R1, R3 and R6.

(\romannumeral5) Suppose $n_4(v)\leq1$. If $n_4(v)=1$, then $ch_1(v)\geq3-3\times\frac{2}{3}-\max\{\frac{4}{9}+\frac{1}{2}, 2\times\frac{4}{9}\}=\frac{1}{18}>0$ by R1 and R3. If $n_4(v)=0$, then $ch_1(v)\geq3-3\times\frac{2}{3}-2\times\frac{4}{9}=\frac{1}{9}>0$ by R1 and R3.
\end{proof}

\begin{claim}\label{W22}
For each vertex $v\in W_2$ with $f_3(v)=2$ and $f_4(v)=0$, $ch_2(v)\ge0$.
\end{claim}
\begin{proof}
Firstly, suppose that the two $3$-faces are consecutive and denote them by $f_1$ and $f_2$. Assume that there exists one $6^+$-face in $\{f_4,f_5,f_6\}$, then $ch_1(v)\geq 3-4\times\frac{2}{3}-\frac{1}{3}=0$ by R1 and R3. Next we consider the situation where $d(f_i)=5$ for each $i\in\{4,5,6\}$, see Figure \ref{fv3eps}$(C_2)$.

If $d(v_i)\geq5$ for some $i\in\{4,5\}$, then $\max\{\send{v}{f_{i-1}}, \send{v}{f_{i}}\}\leq \frac{1}{2}$ by R3, and thus $ch_1(v)\geq3-3\times\frac{2}{3}-2\times\frac{1}{2}=0$ by R1 and R3. Now let $d(v_4)=d(v_5)=4$. Since $S_{3}\not\subseteq G$, $n_{5^+}(f_i)\geq2$ for some $i\in\{3,4\}$ and $n_{5^+}(f_j)\geq2$ for some $j\in\{4,5\}$. If $i\neq j$, then $ch_1(v)\geq3-3\times\frac{2}{3}-2\times\frac{1}{2}=0$ by R1 and R3. If $i=j=4$, then we may assume that $n_{5^+}(f_3)=n_{5^+}(f_5)=1$ (otherwise $ch_1(v)\geq0$). Note that $f_3(v_k)\leq1$ for each $k\in\{4,5\}$. If $v_k$ is good for some $k\in\{4,5\}$, then $ch_1(v_k)\geq\frac{1}{3}$ by R1-R2. Hence, $v_k$ sends at least $\frac{1}{3}$ to $v$ (if $ch_1(v)<0$) via a nice path by R6. Since $S_2,S_{12}\nsubseteq G$, we get that at least one vertex in $\{v_4,v_5\}$ is good, and we are done.

Secondly, suppose that the two $3$-faces are not consecutive, say $f_1$ and $f_3$ are the $3$-faces. By $S_{16}\nsubseteq G$, $\zeta_v(f_{3b})\leq1$. If $d(f_2)\geq6$, then according to $S_3$, we have that $ch_1(v)\geq3-3\times\frac{2}{3}-\frac{1}{2}-\frac{1}{3}-\frac{1}{9}=\frac{1}{18}>0$ by R1, R3 and R5. If $d(f_4)\geq6$ and $\zeta_v(f_{3b})=1$, then $n_{5^+}(f_2)\geq2$ by $S_2\nsubseteq G$, and so $\send{v}{f_2}\leq \frac{1}{2}$ by R3. Thus $ch_1(v)\geq3-3\times\frac{2}{3}-\frac{1}{2}-\frac{1}{3}-\frac{1}{9}=\frac{1}{18}>0$ by R1, R3 and R5. In the following, we may assume $d(f_i)=5$ for each $i\in\{2,4,5\}$, see Figure \ref{fv3eps}($C_3$).

Assume $\zeta_v(f_{3b})=1$, and let $v_1v_2$ be the edge incident with a bad $5$-face. By $S_2\nsubseteq G$, we get $n_{5^+}(f_2)\geq2$ and $n_{5^+}(f_5)\geq2$, and so $\send{v}{f_i}\leq \frac{1}{2}$ for each $i\in\{2,5\}$ by R3. If $f_3(v_i)\leq1$ for some $i\in[2]$, then $v$ need not send any charge to the bad $5$-face by R5 (since $v_i$ sends $\frac{2}{3}$ to the bad $5$-face), and thus $ch_1(v)\geq3-3\times\frac{2}{3}-2\times\frac{1}{2}=0$ by R1 and R3. It remains to consider $f_3(v_i)=2$ for each $i\in[2]$. If $n_{5^+}(f_4)\geq2$, then $\send{v}{f_4}\leq \frac{1}{2}$ by R3, and thus $ch_1(v)\geq3-2\times\frac{2}{3}-3\times\frac{1}{2}-\frac{1}{9}=\frac{1}{18}>0$ by R1, R3 and R5. Otherwise, $n_{5^+}(f_4)=1$. We have that $n_{5^+}(f_2)\geq3$ since $S_2\nsubseteq G$, and $n_{6^{+}}(f_5)\geq1$ or $n_{5^{+}}(f_5)\geq3$ since $S_{28}\nsubseteq G$, and so $\send{v}{f_i}\leq \frac{4}{9}$ for each $i\in\{2,5\}$. Hence, $ch_1(v)\geq 3-3\times \frac{2}{3}-2\times \frac{4}{9}-\frac{1}{9}=0$ by R1, R3 and R5.

Assume $\zeta_v(f_{3b})=0$. Since $S_{3}\not\subseteq G$, we know that at least one of $f_i\in \{f_4,f_5\}$ satisfies $n_{5^+}(f_i)\geq2$. If $n_{5^+}(f_i)\geq2$ for each $i\in\{4,5\}$, then $\send{v}{f_i}\leq \frac{1}{2}$ for each $i\in\{4,5\}$ by R3, and thus $ch_1(v)\geq 3-3\times\frac{2}{3}-2\times\frac{1}{2}=0$ by R1 and R3. Otherwise we assume that $n_{5^+}(f_4)=1$ (which means $n_{5^+}(f_5)\geq2$), then $n_{5^+}(f_2)\geq2$ by $S_2\nsubseteq G$, and thus $ch_1(v)\geq3-3\times\frac{2}{3}-2\times\frac{1}{2}=0$ by R1 and R3.
\end{proof}
\begin{claim}\label{W23}
For each vertex $v\in W_2$ with $f_3(v)=2$ and $f_4(v)=1$, $ch_2(v)\ge0$.
\end{claim}
\begin{proof}
There are two subcases to be considered, see Figure \ref{fv34eps}. Recall that $v$ sends no charge to any bad $5$-face by R5.

\begin{figure}[H]
 \begin{center}
   \includegraphics[scale=0.7]{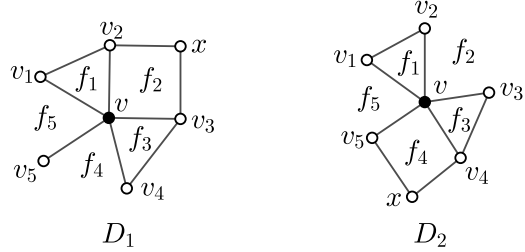}\\
   \caption[14cm]{Configuration for $5$-vertex $v$.}
   \label{fv34eps}
 \end{center}
\end{figure}

We consider the configuration $D_1$ first. (\romannumeral1) Suppose $n_4(v)=5$, that is $d(v_i)=4$ for each $i\in[5]$. Since $S_1 \nsubseteq G$, we obtain that $d(x)\geq5$. If $d(f_i)=5$ for $i\in\{4,5\}$, then by $S_2\nsubseteq G$, $n_{5^+}(f_i)\geq2$. Thus $ch_1(v)\geq 3-3\times\frac{2}{3}-\max\{2\times\frac{1}{2},\frac{1}{2}+\frac{1}{3}\}=0$ by R1 and R3.

(\romannumeral2) Suppose $n_4(v)=4$. If $d(v_1)\geq5$, then $d(x)\geq5$ by $S_1 \nsubseteq G$, and $f_4\neq (4,4,4,4,5)$ by $S_2 \nsubseteq G$. If $d(v_5)\geq5$, then $d(x)\geq5$ by $S_1 \nsubseteq G$. In both cases, $ch_1(v)\geq 3-3\times\frac{2}{3}-\max\{2\times\frac{1}{2},\frac{1}{2}+\frac{1}{3}\}=0$ by R1 and R3. At last, we study the case where $d(v_2)\geq5$. If $f_{6^+}(v)\geq1$, then $ch_1(v)\geq3-4\times\frac{2}{3}-\frac{1}{3}=0$ by R1 and R3. We now turn to the situation $f_{6^+}(v)=0$. In this situation, we may assume that $n_{5^+}(f_5)=1$ (otherwise $ch_1(v)\geq0$). Let us see $v_5$. Note that $f_3(v_5)\leq1$. Denote by $f_6$ and $f_7$ the remaining faces incident with $v_5$ in clockwise. If $v_5$ is good, then $v_5$ sends at least $\frac{1}{3}$ to $v$ (if $ch_1(v)<0$) via a nice path by R6. Otherwise $d(f_6)=3$ and $f_7$ is a bad $5$-face, then by $S_2\nsubseteq G$, we have $d(x_1)\geq5$, see Figure \ref{C4eps}$(E_1)$. By the assumption, $d(v_2)\geq5$. If $d(v_2)\geq6$, then $v$ sends at most $\frac{5}{9}$ to $f_2$ by R4; if $d(v_2)=5$, then by $S_{6}\nsubseteq G$, $d(x)\geq5$, and $v$ sends at most $\frac{5}{9}$ to $f_2$ by R3. On the other hand, if $d(x_1)\geq6$, then $v$ sends at most $\frac{4}{9}$ to $f_4$ by R3; if $d(x_1)=5$, then by $S_{7}\nsubseteq G$, $d(x_2)\geq5$, and $v$ sends at most $\frac{4}{9}$ to $f_4$ by R3. In conclusion, $ch_1(v)\geq3-3\times\frac{2}{3}-\frac{4}{9}-\frac{5}{9}=0$ by R1 and R3.

(\romannumeral3) Suppose $n_4(v)=3$, that is $n_{5^+}(v)=2$. If the pair of two $5^+$-vertices fall in $\{(v_3,v_4),(v_3,v_5),(v_4,v_5)\}$, then by $S_1,S_2\nsubseteq G$, we get $ch_1(v)\geq3-3\times\frac{2}{3}-\max\{\frac{1}{2}+\frac{4}{9}, 2\times\frac{1}{2}\}=0$ by R1 and R3. It remains to consider the pairs: $\{(v_1,v_4)$, $(v_2,v_4)$, $(v_2,v_3)\}$.

Assume $d(v_2)\geq5$ and $d(v_4)\geq5$. If $f_{6^+}(v)\geq1$, then $ch_1(v)\geq3-3\times\frac{2}{3}-\max\{\frac{2}{3}+\frac{1}{3}, \frac{1}{2}+\frac{1}{3}\}=0$ by R1 and R3. It remains to discuss the case where $f_{6^+}(v)=0$. Here, we can let $n_{5^+}(f_5)=1$ (otherwise $ch_1(v)\geq0$) and $v_5$ be not good (otherwise $v_5$ could send at least $\frac{1}{3}$ to $v$ (if $ch_1(v)<0$) via a nice path by R6 and $ch_2(v)\geq0$). Denote by $f_6$ and $f_7$ the remaining faces incident with $v_5$ in clockwise. Note that $d(f_6)=3$ and $f_7$ is a bad $5$-face. By $S_2\nsubseteq G$, we get $n_{5^+}(f_4)\geq3$, and so $\send{v}{f_4}\leq \frac{4}{9}$. On the other hand, recall that $d(v_2)\geq5$. If $d(v_2)\geq6$, then $v$ sends at most $\frac{5}{9}$ to $f_2$ by R3; if $d(v_2)=5$, then $d(x)\geq5$ holds because of $S_{21}\nsubseteq G$. Hence, $ch_1(v)\geq3-3\times\frac{2}{3}-\frac{4}{9}-\frac{5}{9}=0$ by R1 and R3.

Assume $d(v_2)\geq5$ and $d(v_3)\geq5$. If $f_{6^{+}}(v)=1$, then $ch_1(v)\geq3-3\times\frac{2}{3}-\frac{5}{9}-\frac{1}{3}=\frac{1}{9}>0$ by R1 and R3. Suppose that $f_{6^{+}}(v)=0$. Since $S_3\nsubseteq G$, $n_{5^+}(f_i)\geq2$ holds for some $i\in \{4,5\}$. If $n_{5^+}(f_i)\geq2$ for each $i\in \{4,5\}$, then $\send{v}{f_i}\leq \frac{1}{2}$ for each $i\in\{4,5\}$ by R3, and thus $ch_1(v)\geq3-2\times\frac{2}{3}-2\times\frac{1}{2}-\frac{5}{9}=\frac{1}{9}>0$ by R1 and R3. Otherwise let $n_{5^+}(f_5)=1$, that is $d(x_1)=d(x_2)=4$, see Figure \ref{C4eps}$(E_2)$. If $d(x)\geq5$, then $\send{v}{f_2}\leq \frac{2}{3}$ by R3, and thus  $ch_1(v)\geq3-3\times\frac{2}{3}-2\times\frac{1}{2}=0$ by R1 and R3. Otherwise if $d(x)=4$, then $f_3(v_5)\leq1$ as $f_2$ is a $4$-face and any two $4$-faces in $G$ are at distance at least $2$. If $v_5$ is good, then $v_5$ sends at least $\frac{1}{3}$ to $v$ (if $ch_1(v)<0$) via a nice path by R6. Otherwise $f_3(v_5)=1$ and $v_5$ is incident with a bad $5$-face. Denote by $f_6$ and $f_7$ the faces incident with $v_5$ in clockwise. Since $S_2\nsubseteq G$, we get $d(f_6)=3$ and $f_7$ is a bad $5$-face. If $d(z_2)=5$, then $n_{5^+}(f_4)\geq3$ by $S_{14} \nsubseteq G$, and so $\send{v}{f_4}\leq \frac{4}{9}$. Otherwise $d(z_2)\geq6$, in this situation $\send{v}{f_4}\leq \frac{4}{9}$ by R3. Hence,  $ch_1(v)\geq3-3\times\frac{2}{3}-\frac{5}{9}-\frac{4}{9}=0$ by R1 and R3.
\begin{figure}[H]
 \begin{center}
   \includegraphics[scale=0.7]{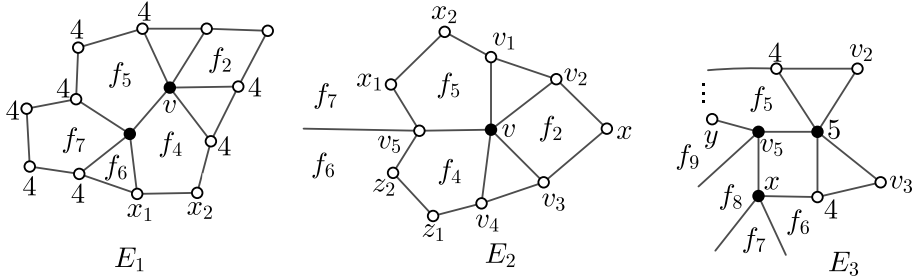}\\
   \caption[14cm]{Specified Configuration.}
   \label{C4eps}
 \end{center}
\end{figure}

Finally, we consider the case where $d(v_1)\geq5$ and $d(v_4)\geq5$.
$(a)$. $d(x)=4$. Let us start to claim that $v$ sends at most $\frac{8}{9}$ in total to $\{f_4,f_5\}$. Assume that $f_{6^+}(v)=0$. If $d(v_1)\geq6$, then we are done by R3. Otherwise if $d(v_1)=5$, then $n_{5^+}(f_5)\geq3$ holds because of $S_{8}\nsubseteq G$, and so $\send{v}{f_5}\leq \frac{4}{9}$. The above arguments can also be applied to $v_4$. So the same result holds for $f_{6^+}(v)\geq1$, as claimed. Note that $v_1$ is symmetric to $v_4$. So we only discuss $v_4$ in the following, and we would like to claim that $v_4$ could send at least $\frac{1}{9}$ to $v$ (if $ch_1(v)<0$) when $d(v_4)\geq6$ via a nice path.

Assume $d(v_4)\geq7$. By Claim \ref{c3}, we have $ch_1(v_4)\geq\frac{5d(v)-36}{18}+\frac{1}{9}+\frac{2}{9}=\frac{5d(v)-30}{18}>\frac{1}{9}$. Assume $d(v_4)=6$. Since $S_{45}\not\subseteq G$, $\zeta_{v_4}(f_{3b})\leq1$. If $\zeta_{v_4}(f_{3b})=0$, then $ch_1(v_4)\geq\frac{5d(v)-36}{18}+4\times\frac{1}{9}\geq\frac{1}{9}$. Otherwise if $\zeta_{v_4}(f_{3b})=1$, then by $S_{43}\nsubseteq G$, we get $n_{5^+}(f_4)\geq3$ and thus $v_4$ sends at most $\frac{4}{9}$ to $f_4$ by R4. Hence $ch_1(v_4)\geq\frac{5d(v)-36}{18}+2\times\frac{1}{9}+\frac{2}{9}\geq\frac{1}{9}$, and $v_4$ could send at least $\frac{1}{9}$ to $v$ (if $ch_1(v)<0$) via a nice path by R6, as claimed. So when $\min\{d(v_2),d(v_4)\}\geq6$, $v$ (if $ch_1(v)<0$) could receive at least $\frac{2}{9}$ in total from $\{v_1,v_4\}$ via two nice paths by R6, and $ch_2(v)\geq3-2\times\frac{2}{3}-1-2\times\frac{4}{9}+2\times\frac{1}{9}=0$ by R1, R3 and R6.
\begin{figure}[H]
 \begin{center}
   \includegraphics[scale=0.7]{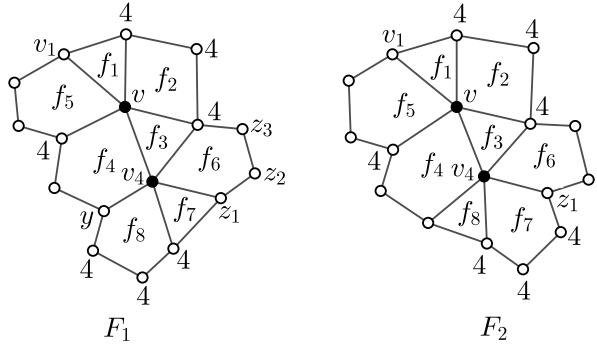}\\
   \caption[14cm]{Specified Configuration.}
   \label{5-vx1eps}
 \end{center}
\end{figure}
Now we consider $\min\{d(v_2),d(v_4)\}=5$. W.l.o.g., we assume $d(v_4)=5$. Since $S_{17} \nsubseteq G$, $\zeta_{v_4}(f_{3b})=0$. Denote by $f_6$, $f_7$ and $f_8$ the remaining faces incident with $v_4$ in clockwise. As $S_{10}\nsubseteq G$, $n_{5^+}(f_6)\geq2$ when $d(f_6)=5$. We may assume that $f_{6^+}(v_4)=0$ (otherwise $ch_1(v_4)\geq3-\frac{2}{3}-\frac{1}{3}-\max\{\frac{2}{3}+2\times\frac{1}{2}, \frac{2}{3}+\frac{1}{2}+\frac{4}{9}\}>\frac{2}{9}$ by R1 and R3, and $v_4$ could send at least $\frac{2}{9}$ to $v$ (if $ch_1(v)<0$) via a nice path by R6. So $ch_2(v)\geq0$). First, let $d(f_7)=3$. If $n_{5^+}(f_8)=1$, see Figure \ref{5-vx1eps}$(F_1)$ ($d(y)=4$), then by $S_{8},S_{13}\nsubseteq G$, $y$ can not be incident with a bad $5$-face. Note that $f_3(y)\leq1$, thus $y$ is good, and $ch_1(y)\geq\frac{1}{3}$. Hence, $y$ could send at least $\frac{1}{3}$ to $v$ (if $ch_1(v)<0$) via a nice path by R6. Otherwise if $n_{5^+}(f_8)\geq2$, then  $ch_1(v_4)\geq3-2\times\frac{2}{3}-2\times\frac{1}{2}-\frac{4}{9}=\frac{2}{9}$ by R1 and R3. Second, let $d(f_8)=3$. If $n_{5^+}(f_7)=1$, see Figure \ref{5-vx1eps}$(F_2)$ ($d(z_1)=4$), then $z_1$ is good since $S_{13} \nsubseteq G$, and thus $z_1$ could send at least $\frac{1}{3}$ to $v$ (if $ch_1(v)<0$) via a nice path by R6. Otherwise if $n_{5^+}(f_7)\geq2$, then $ch_1(v_4)\geq3-2\times\frac{2}{3}-2\times\frac{1}{2}-\frac{4}{9}=\frac{2}{9}$ by R1 and R3.

In conclusion, when $d(v_4)=5$, $v$ (if $ch_1(v)<0$) could receive at least $\frac{2}{9}$ from one vertex in $\{v_4,y_1,z_1\}$ via a nice path by R6. Thus $ch_2(v)\geq3-2\times\frac{2}{3}-1-2\times\frac{4}{9}+\frac{2}{9}=0$ by R1, R3 and R6.

$(b)$. $d(x)\geq5$. Then $ch_1(v)\geq 3-3\times\frac{2}{3}-2\times\frac{1}{2}=0$ by R1 and R3.

(\romannumeral4) Suppose $n_4(v)=2$. If the pair of two $4$-vertices fall in  $\{(v_1,v_2),(v_1,v_3)$, $(v_1,v_4),$ $(v_2,v_5)\}$, then $ch_1(v)\geq3-3\times\frac{2}{3}-2\times\frac{1}{2}=0$ by R1 and R3.

Assume $d(v_1)=d(v_5)=4$. If $f_{6^+}(v)\geq1$, then $ch_1(v)\geq3-3\times\frac{2}{3}-\frac{5}{9}-\frac{1}{9}=\frac{1}{3}>0$ by R1 and R3. It remains to consider $f_{6^+}(v)=0$. If $n_{5^+}(f_5)\geq2$, then $\send{v}{f_5}\leq \frac{1}{2}$ by R3, and thus  $ch_1(v)\geq3-2\times\frac{2}{3}-\frac{5}{9}-2\times\frac{1}{2}=\frac{1}{9}>0$ by R1 and R3. Otherwise $n_{5^+}(f_5)=1$. We may let $v_5$ is not good (otherwise $v_5$ could send $\frac{1}{3}$ to $v$ (if $ch_1(v)<0$) via a nice path by R6, and thus $ch_2(v)\geq3-3\times\frac{2}{3}-\frac{5}{9}-\frac{1}{2}+\frac{1}{3}>0$ by R1, R3 and R6). Denote by $f_6$ and $f_7$ the faces incident with $v_5$ in clockwise. Since $S_2 \nsubseteq G$, $d(f_6)=3$ and $f_7$ is a bad $5$-face. Moreover, by $S_2\nsubseteq G$ again, $n_{5^+}(f_4)\geq3$, and so $\send{v}{f_4}\leq \frac{4}{9}$.  Thus $ch_1(v)\geq3-3\times\frac{2}{3}-\frac{5}{9}-\frac{4}{9}=0$ by R1 and R3.

Assume $d(v_2)=d(v_3)=4$.
$(a)$. $d(x)=4$. By the same arguments as the case $d(v_1)\geq5$ and $d(v_4)\geq5$, we have that $v$ (if $ch_1(v)<0$) could receive at least $\frac{2}{9}$ from $\{v_2,x\}$ via a nice path, and thus $ch_2(v)\geq3-1-2\times\frac{2}{3}-2\times\frac{4}{9}+\frac{2}{9}=0$ by R1, R3 and R6.
$(b)$. $d(x)\geq5$. Then $ch_1(v)\geq3-3\times\frac{2}{3}-2\times\frac{4}{9}=\frac{1}{9}>0$ by R1 and R3.

(\romannumeral5) Suppose $n_4(v)\leq1$. If $n_4(v)=1$, then $ch_1(v)\geq3-2\times\frac{2}{3}-\max\{\frac{1}{2}+\frac{4}{9}+\frac{5}{9}, \frac{2}{3}+2\times\frac{4}{9}, 2\times\frac{1}{2}+\frac{5}{9}\}=\frac{1}{9}>0$ by R1 and R3. If $n_4(v)=0$, then $ch_1(v)\geq3-2\times\frac{2}{3}-2\times\frac{4}{9}-\frac{5}{9}=\frac{2}{9}>0$ by R1 and R3.

Now we consider the configuration $D_2$. (\romannumeral1) Suppose $n_4(v)=5$. Since $S_2 \nsubseteq G$, $f_i\neq (4,4,4,4,5)$ for each $i\in\{2,5\}$. By $S_1\nsubseteq G$, we get $f_4\neq(4,4,4,4,5)$, and so $\send{v}{f_4}\leq \frac{1}{2}$ by R3. Thus $ch_1(v)\geq3-3\times\frac{2}{3}-\max\{2\times\frac{1}{2},\frac{1}{2}+\frac{1}{3}\}=0$ by R1 and R3.

(\romannumeral2) Suppose $n_4(v)=4$, that is $n_{5^+}(v)=1$. Assume $d(v_1)\geq5$, then $f_4\neq (4,4,4,4,5)$ holds because of $S_1\nsubseteq G$. Moreover, $f_2\neq(4,4,4,4,5)$. Thus $ch_1(v)\geq3-3\times\frac{2}{3}-2\times\frac{1}{2}=0$ by R1 and R3.

Assume $d(v_2)\geq5$, then $f_4\neq(4,4,4,4,5)$ holds by $S_1\nsubseteq G$, and so $\send{v}{f_4}\leq \frac{1}{2}$ by R3. If $f_{6^{+}}(v)\geq1$, then $ch_1(v)\geq 3-4\times\frac{2}{3}-\frac{1}{3}=0$ by R1 and R3. Now we discuss $d(f_2)=d(f_5)=5$. We may assume that $n_{5^+}(f_5)=1$ (otherwise $ch_1(v)\geq3-3\times\frac{2}{3}-2\times\frac{1}{2}=0$). On the other hand, we may assume that $v_5$ is not good (otherwise $v_5$ could send at least $\frac{1}{3}$ to $v$ (if $ch_1(v)<0$) via a nice path and thus $ch_2(v)\geq3-4\times\frac{2}{3}-\frac{1}{2}+\frac{1}{3}>0$). Denote by $f_6$ and $f_7$ the faces incident with $v_5$ in clockwise. Since $S_2 \nsubseteq G$, we have that $d(f_6)=3$ and $f_7$ is a bad $5$-face. By $S_{18}\nsubseteq G$, we know that $n_{6^+}(f_4)=1$, and thus $v$ sends at most $\frac{5}{9}$ to $f_4$ by R3. Next we claim that $v$ sends at most $\frac{4}{9}$ to $v$. If $d(v_2)\geq6$, then we are done by R3; if $d(v_2)=5$, then $n_{5^+}(f_2)\geq3$ by $S_{11}\nsubseteq G$, and so $\send{v}{f_2}\leq \frac{4}{9}$ by R3, as claimed. Hence, $ch_1(v)\geq3-3\times\frac{2}{3}-\frac{4}{9}-\frac{5}{9}=0$ by R1 and R3.

Assume $d(v_3)\geq5$.
$(a)$. $d(x)=4$. Since $S_2 \nsubseteq G$, $f_5\neq(4,4,4,4,5)$. For brevity, denote by $f_6$ and $f_7$ the faces incident with $v_5$ in clockwise. We may assume that $v_5$ is not good (otherwise, $v_5$ sends at least $\frac{1}{3}$ to $v$ (if $ch_1(v)<0$) via a nice path, and $ch_2(v)\geq3-2\times\frac{2}{3}-1-2\times\frac{1}{2}+\frac{1}{3}=0$). Since $S_2\nsubseteq G$, $d(f_7)=3$ and $f_6$ is a bad $5$-face, which is impossible because $S_{5}$ is reducible.
$(b)$. $d(x)\geq5$. By $S_{2}\nsubseteq G$, we get $f_5\neq(4,4,4,4,5)$, so $ch_1(v)\geq3-3\times\frac{2}{3}-2\times\frac{1}{2}=0$ by R1 and R3.

Assume $d(v_4)\geq5$. If $f_{6^+}(v)\geq1$, then $ch_1(v)\geq3-4\times\frac{2}{3}-\frac{1}{3}=0$ by R1 and R3. Now we consider $f_{6^+}(v)=0$. Notice that $n_{5^+}(f_2)\geq2$ and $n_{5^+}(f_5)\geq2$ holds because of $S_2\nsubseteq G$, and so $\send{v}{f_i}\leq \frac{1}{2}$ for each $i\in\{2,5\}$ by R3. Thus $ch_1(v)\geq3-3\times\frac{2}{3}-2\times\frac{1}{2}=0$ by R1 and R3. Assume $d(v_5)\geq5$, then by $S_2\nsubseteq G$, we get $n_{5^+}(f_2)\geq2$, and so $\send{v}{f_2}\leq \frac{1}{2}$ by R3. Thus $ch_1(v)\geq3-3\times\frac{2}{3}-2\times\frac{1}{2}=0$.

(\romannumeral3) Suppose $n_4(v)=3$, that is $n_{5^+}(v)=2$. If the pair of two $5^+$-vertices fall in $\{(v_1,v_2),(v_2,v_5),(v_3,v_4),(v_3,v_5),(v_4,v_5)\}$, then by $S_2\nsubseteq G$, $ch_1(v)\geq 3-2\times\frac{2}{3}-\max\{2\times\frac{1}{2}+\frac{5}{9}, 2\times\frac{1}{2}+\frac{2}{3}\}=0$ by R1 and R3.

Assume $d(v_2)\geq5$ and $d(v_4)\geq5$. If $f_{6^+}(v)\geq1$, then $ch_1(v)\geq3-4\times\frac{2}{3}-\frac{1}{3}=0$ by R1 and R3. Next we discuss $f_{6^+}(v)=0$. We may assume that $n_{5^+}(f_5)=1$ (otherwise $ch_1(v)\geq0$).
$(a)$. $d(x)=4$. By $S_9\nsubseteq G$, we get $d(v_2)\geq6$, and then $\send{v}{f_2}\leq \frac{4}{9}$ by R3. By $S_{15}\nsubseteq G$, we know that $d(v_4)\geq6$, and then $\send{v}{f_4}\leq \frac{5}{9}$ by R3. Thus $ch_1(v)\geq3-3\times\frac{2}{3}-\frac{5}{9}-\frac{4}{9}=0$ by R1 and R3.
$(b)$. $d(x)\geq5$. Similarly as above, we have $ch_1(v)\geq3-3\times\frac{2}{3}-\frac{5}{9}-\frac{4}{9}=0$ by R1 and R3.

Assume $d(v_2)\geq5$ and $d(v_3)\geq5$.
$(a)$. $d(x)=4$. Let $y$ be the neighbor of $v_5$ which locates on $f_5$. If $d(y)=4$, then by $S_2\nsubseteq G$, $v_5$ is good and thus $v_5$ could send at least $\frac{1}{3}$ to $v$ (if $ch_1(v)<0$) via a nice path by R6. When $d(f_5)\geq6$, and we have $ch_2(v)\geq3-1-2\times\frac{2}{3}-\frac{1}{3}-\frac{4}{9}+\frac{1}{3}=\frac{2}{9}>0$ by R1 and R3. When $d(f_5)=5$, $n_{5^+}(f_5)\geq2$ holds because of $S_4\nsubseteq G$, and we have $\send{v}{f_5}\leq \frac{1}{2}$ by R3. Thus  $ch_2(v)\geq3-1-2\times\frac{2}{3}-\frac{1}{2}-\frac{4}{9}+\frac{1}{3}=\frac{1}{18}>0$ by R1 and R3. It remains to consider $d(y)\geq5$. Note that $f_5\neq(4,4,4,4,5)$. Denote by $f_8$ and $f_9$ the faces incident with $v_5$ in clockwise, and $f_6$, $f_7$ the remaining faces incident with $x$ in clockwise, see Figure \ref{C4eps}$(E_3)$. We may assume that $v_5$ is not good (otherwise $v_5$ could send at least $\frac{1}{3}$ to $v$ via a nice path and $ch_2(v)\geq0$). We immediately have $d(f_9)=3$ and $f_8$ is a bad $5$-face. In this situation, $f_3(x)\leq1$, and by $S_2\nsubseteq G$, $x$ is good. Thus $x$ could send at least $\frac{1}{3}$ to $v$ (if $ch_1(v)<0$) via a nice path and $ch_2(v)\geq0$.

$(b)$. $d(x)\geq5$. If $f_{6^+}(v)\geq1$, then $ch_2(v)\geq3-4\times\frac{2}{3}-\frac{1}{3}=0$ by R1 and R3.  It remains to discuss $f_{6^+}(v)=0$. In this situation, we may assume that $d(y)=4$ (otherwise $ch_2(v)\geq3-3\times\frac{2}{3}-\frac{1}{2}-\frac{4}{9}>0$). Denote by $f_6$ and $f_7$ the remaining faces incident with $v_5$ in clockwise. Let $v_5$ be a vertex which is not good (otherwise $v_5$ could send at least $\frac{1}{3}$ to $v$ (if $ch_1(v)<0$) via a nice path and $ch_2(v)\geq3-4\times\frac{2}{3}-\frac{4}{9}+\frac{1}{3}>0$). Note that $d(f_6)=3$ and $f_7$ is a bad $5$-face by $S_2\nsubseteq G$. Since $S_{18}\nsubseteq G$, we have that $d(x)\geq6$, and $v$ sends at moat $\frac{5}{9}$ to $f_4$ by R3. Hence, $ch_2(v)\geq3-3\times\frac{2}{3}-\frac{5}{9}-\frac{4}{9}=0$ by R1 and R3.

Assume $d(v_1)\geq5$ and $d(v_3)\geq5$. $(a)$. $d(x)=4$. Denote by $f_6$, $f_7$ and $f_8$ the remaining faces incident with $x$ in clockwise and $f_9$ another faces incident with $v_5$. If $v_5$ is good, then $v$ (if $ch_1(v)<0$) could receive at least $\frac{1}{3}$ from $v_5$ via a nice path by R6. Now we assume $v_5$ is not good. Since $S_2\nsubseteq G$, $d(f_9)=3$ and $f_8$ is a bad $5$-face. In this situation, $x$ must be good by $S_2\nsubseteq G$ again. Hence, $v$ (if $ch_1(v)<0$) could receive at least $\frac{1}{3}$ from $x$ via a nice path by R6. Thus $ch_2(v)\geq3-1-2\times\frac{2}{3}-2\times\frac{1}{2}+\frac{1}{3}=0$ by R1, R3 and R6.
$(b)$. $d(x)\geq5$. Then $ch_2(v)\geq3-3\times\frac{2}{3}-2\times\frac{1}{2}=0$ by R1 and R3.

(\romannumeral4) Suppose $n_4(v)=2$. If the pair of two $4$-vertices fall in $\{(v_1,v_2),(v_1,v_3),(v_1,v_4),(v_2,v_3),$ $(v_2,v_4),(v_2,v_5),(v_3,v_5)\}$, then $ch_2(v)\geq3-2\times\frac{2}{3}-\max\{2\times\frac{1}{2}+\frac{5}{9}, \frac{2}{3}+\frac{1}{2}+\frac{4}{9},\frac{2}{3}+\frac{5}{9}+\frac{4}{9},\frac{2}{3}+2\times\frac{1}{2}\} =0$ by R1 and R3.

Assume $d(v_1)=d(v_5)=4$. If $f_{6^+}(v)\geq1$, then $ch_2(v)\geq3-4\times\frac{2}{3}-\frac{1}{3}=0$ by R1 and R3. Otherwise $f_{6^+}(v)=0$, then we may assume that $n_{5^+}(f_5)=1$ (otherwise $ch_2(v)\geq0$). If $d(x)\geq5$ or $d(v_4)\geq6$, then $\send{v}{f_4}\leq \frac{4}{9}$ by R3, and thus $ch_2(v)\geq3-3\times\frac{2}{3}-\frac{4}{9}-\frac{5}{9}=0$ by R1 and R3. So we consider $d(x)=4$ and $d(v_4)=5$. Denote by $f_6$ and $f_7$ the remaining faces incident with $v_5$ in clockwise. Since $S_2\nsubseteq G$, $v_5$ is good. Hence, $v_5$ could send at least $\frac{1}{3}$ to $v$ (if $ch_1(v)<0$) via a nice path by R6, and $ch_2(v)\geq3-4\times\frac{2}{3}-\frac{4}{9}+\frac{1}{3}>0$ by R1, R3 and R6.

Assume $d(v_4)=d(v_5)=4$. $(a)$. $d(x)=4$. Denote by $f_6$, $f_7$ the remaining faces incident with $v_5$ in clockwise. Since $S_2\nsubseteq G$, $d(f_7)=3$ and $f_6$ is a bad $5$-face. We also have $n_{5^+}(f_5)\geq3$ because of $S_2\nsubseteq G$. Denote by $f_8$, $f_9$ the remaining faces incident with $x$ in clockwise. By $S_2\nsubseteq G$, $x$ is good and $x$ could send at least $\frac{1}{3}$ to $v$ (if $ch_1(v)<0$) via a nice path. Thus $ch_2(v)\geq3-1-2\times\frac{2}{3}-2\times\frac{4}{9}+\frac{1}{3}>0$ by R1 and R3.
$(b)$. $d(x)\geq5$. Then $ch_2(v)\geq3-3\times\frac{2}{3}-\frac{1}{2}-\frac{4}{9}=\frac{1}{18}>0$ by R1 and R3.

(\romannumeral5) Suppose $n_4(v)\leq1$. If $n_4(v)=1$, then $ch_2(v)\geq3-2\times\frac{2}{3}-\max\{\frac{1}{2}+\frac{5}{9}+\frac{4}{9},\frac{2}{3}+2\times\frac{4}{9}, \frac{2}{3}+\frac{1}{2}+\frac{4}{9}\}=\frac{1}{18}>0$ by R1 and R3. If $n_4(v)=0$, then $ch_2(v)\geq3-2\times\frac{2}{3}-\frac{5}{9}-2\times\frac{4}{9}=\frac{2}{9}>0$ by R1 and R3.
\end{proof}

\begin{claim}\label{W24}
For each vertex $v\in W_2$ with $f_3(v)=1$, $ch_2(v)\ge0$.
\end{claim}
\begin{proof}
 W.l.o.g., let $d(f_1)=3$.

\textbf{Case 1.} Suppose $f_4(v)=0$. Assume that $\zeta_{v}(f_{3b})=0$ firstly. If $f_{6^{+}}(v)=1$, then $ch_1(v)\geq3-4\times\frac{2}{3}-\frac{1}{3}=0$ by R1 and R3. Otherwise $f_{6^{+}}(v)=0$. Since $S_3\nsubseteq G$, there exists at least one $5$-face $f_i$ ($i\in\{2,3,4,5\}$) such that $n_{5^+}(f_i)\geq2$, and $\send{v}{f_i}\leq \frac{1}{2}$ by R3. If $n_{5^+}(f_1)=3$, then $\send{v}{f_i}\leq \frac{1}{2}$ for each $i\in\{2,5\}$ by R3, and thus $ch_1(v)\geq3-2\times\frac{2}{3}-3\times\frac{1}{2}=\frac{1}{6}>0$ by R1 and R3. If $n_{5^+}(f_1)=2$, say $d(v_1)\geq5$, then we claim that $v$ sends at most $\frac{10}{9}$ in total to $\{f_1, f_5\}$. Obviously, the claim holds when $d(v_1)\geq6$; when $d(v_1)=5$, then by $S_{11}\nsubseteq G$, either $n_{5^+}(f_1)\geq2$ or $n_{5^+}(f_5)\geq3$, as claimed. Hence, $ch_1(v)\geq3-2\times\frac{2}{3}-\frac{1}{2}-\frac{10}{9}=\frac{1}{18}>0$ by R1 and R3. If $n_{5^+}(f_1)=0$, that is $d(v_1)=d(v_2)=4$, then by $S_2\nsubseteq G$, $n_{5^+}(f_1)\geq2$ and $n_{5^+}(f_5)\geq2$, and so $\send{v}{f_i}\leq \frac{1}{2}$ for each $i\in\{2,5\}$ by R3. Thus $ch_1(v)\geq3-2\times\frac{2}{3}-3\times\frac{1}{2}=\frac{1}{6}>0$ by R1 and R3.

We now turn to the case $\zeta_{v}(f_{3b})=1$, which means that $d(v_1)=d(v_2)=4$. Since $S_2,S_3\nsubseteq G$, we get that $ch_1(v)\geq3-\frac{2}{3}-\frac{1}{2}-\max\{\frac{1}{3}+\frac{1}{2}+\frac{2}{3}, 2\times\frac{1}{2}+\frac{2}{3}\}-\frac{1}{9}=\frac{1}{18}>0$ by R1, R3 and R5.

\textbf{Case 2.} Suppose $f_4(v)=1$ and let the other vertex on $4$-face is $x$. Recall that $v$ sends no charge to a bad $5$-face (if it exists) which is incident with a $(4,4,v)$-face by R5. By symmetry, we only need to consider the cases $d(f_2)=4$ and $d(f_3)=4$.

\textbf{Subcase 2.1.} $f_{6^+}(v)=1$.
$(a)$. $n_{5^+}(f_i)=1$ for some $i\in\{2,3\}$. Then we have $ch_1(v)\geq 3-1-3\times\frac{2}{3}-\frac{1}{3}=-\frac{1}{3}$. If $ch_1(v)\geq0$, then we are done. So $ch_1(v)<0$, that is, $v$ is poor. Clearly, if there is a good $4$-vertex in $N(v)$, then $ch_2(v)\geq -\frac{1}{3}+\frac{1}{3}=0$ by Claim \ref{c2} and R6. Next we discuss the case that there is no good $4$-vertex in $N(v)$.

Assume that $d(f_2)=4$, then $d(v_1)\geq5$ as $S_1\nsubseteq G$. If $d(f_3)\geq6$, then we may assume that $d(v_5)=4$ (otherwise if $d(v_5)\geq 5$, then $ch_1(v)\geq3-1-\frac{2}{3}-\frac{1}{3}-2\times \frac{1}{2}=0$). Recall that $v_5\in N(v)$ is not good. By $S_2,S_{23}\nsubseteq G$, we get that $n_{5^+}(f_4)\geq2$, and then $\send{v}{f_4}\leq \frac{1}{2}$ by R3. Thus $ch_1(v)\geq3-1-\frac{1}{3}-\frac{2}{3}-2\times\frac{1}{2}=0$ by R1 and R3. If $d(f_4)\geq6$, then $n_{5^+}(f_3)\geq2$ because of $S_{4}\nsubseteq G$, and so $\send{v}{f_3}\leq \frac{1}{2}$ by R3. Thus $ch_1(v)\geq3-1-\frac{1}{3}-\frac{2}{3}-2\times\frac{1}{2}=0$ by R1 and R3. If $d(f_5)\geq6$, then we may assume that $d(v_4)=4$ (otherwise $ch_1(v)\geq0$) and $v_4\in N(v)$ is not good. Similarly, by $S_2,S_{23}\nsubseteq G$, we get that $n_{5^+}(f_4)\geq2$, and then $\send{v}{f_4}\leq \frac{1}{2}$ by R3. Thus $ch_1(v)\geq3-1-\frac{1}{3}-\frac{2}{3}-2\times\frac{1}{2}=0$.

Assume $d(f_3)=4$, then by $S_2,S_3\nsubseteq G$, there is at least one $j\in\{2,4,5\}$ such that $f_j\neq (4,4,4,4,5)$, and so $\send{v}{f_j}\leq \frac{1}{2}$. If $d(f_2)\geq6$, then we may assume that $d(v_5)=4$ (otherwise $ch_1(v)\geq0$) and $v_5\in N(v)$ is not good. Since $S_2,S_{23}\nsubseteq G$, we get that $n_{5^+}(f_5)\geq2$, and $ch_1(v)\geq3-1-\frac{1}{3}-\frac{2}{3}-2\times\frac{1}{2}=0$ by R1 and R3. If $d(f_4)\geq6$, then by the similar arguments, $n_{5^+}(f_5)\geq2$, and thus $ch_1(v)\geq3-1-\frac{1}{3}-\frac{2}{3}-2\times\frac{1}{2}=0$. If $d(f_5)\geq6$, then $n_{5^+}(f_2)\geq2$ and $n_{5^+}(f_4)\geq2$ hold because of $S_{4}\nsubseteq G$, and so $\send{v}{f_i}\leq \frac{1}{2}$ for each $i\in\{2,4\}$ by R3. Thus $ch_1(v)\geq3-1-\frac{1}{3}-\frac{2}{3}-2\times\frac{1}{2}=0$ by R1 and R3.

$(b)$. $n_{5^+}(f_i)\geq2$ for some $i\in\{2,3\}$. Then $ch_1(v)\geq3-4\times\frac{2}{3}-\frac{1}{3}=0$ by R1 and R3.

\textbf{Subcase 2.2.} $f_{6^{+}}(v)=0$.
$(a)$. $n_{5^+}(f_i)=1$ for some $i\in\{2,3\}$. If $d(f_2)=4$, then by $S_1\nsubseteq G$, $d(v_1)\geq5$ holds, and by $S_3\nsubseteq G$, there is at least one face $f_j$ ($j\in \{3,4\}$) satisfying $f_j\neq(4,4,4,4,5)$, and so $\send{v}{f_j}\leq \frac{1}{2}$ and  $\send{v}{f_5}\leq \frac{1}{2}$ by R3. If $d(f_3)=4$, then $f_i\neq (4,4,4,4,5)$ holds for each $i\in\{2,4\}$ by $S_4\nsubseteq G$, and so $\send{v}{f_i}\leq \frac{1}{2}$. By the similar arguments as above, we may assume that each vertex in $N(v)$ is not good (otherwise $ch_2(v)\geq 3-1-2\times \frac{2}{3}-2\times \frac{1}{2}+\frac{1}{3}=0$, and we are done).

Assume that $d(f_2)=4$. For brevity, let $x_1\in N(v_3)$ such that $x_1$ locates on $f_3$, and denote by $f_6$, $f_7$ the faces incident with $v_3$ in clockwise. Since $S_2\nsubseteq G$, $d(f_7)=3$ and $f_6$ is a bad $5$-face, and $d(x_1)\geq5$. Moreover, if $d(x_1)=5$, then $n_{5^+}(f_3)\geq3$ because of $S_{19}\nsubseteq G$. So $v$ sends at most $\frac{4}{9}$ to $f_3$ by R3. However, $x$ is good in this situation and $ch_1(x)\geq\frac{1}{3}$. Hence, $x$ could send at least $\frac{1}{3}$ to $v$ (if $ch_1(v)<0$) via a nice path by R6, and $ch_2(v)\geq3-1-2\times\frac{2}{3}-\frac{1}{2}-\frac{4}{9}+\frac{1}{3}>0$ by R1, R3 and R6.

Assume that $d(f_3)=4$. Note that $v_3, v_4\in N(v)$ are not good. By $S_2\nsubseteq G$, $x$ is good and $ch_1(x)\geq\frac{1}{3}$. Hence, $x$ could send at least $\frac{1}{3}$ to $v$ (if $ch_1(v)<0$) via a nice path by R6, and $ch_2(v)\geq3-1-2\times\frac{2}{3}-2\times\frac{1}{2}+\frac{1}{3}=0$ by R1, R3 and R6.

$(b)$. $n_{5^+}(f_i)\geq2$ for some $i\in\{2,3\}$. Similarly, we may assume that each vertex in $N(v)$ is not good (otherwise $ch_2(v)\geq 3-\frac{5}{9}-4\times \frac{2}{3}+\frac{1}{3}=\frac{1}{9}>0$, and we are done). Assume $d(f_2)=4$. Since $S_3\nsubseteq G$, there exists at least one face $f_i$ and $f_j$ in $\{f_3,f_4\}$ and $\{f_4,f_5\}$, respectively such that $n_{5^+}(f_i)\geq2$, $n_{5^+}(f_j)\geq2$. If $i=j=4$, that is $n_{5^+}(f_3)=n_{5^+}(f_5)=1$, recall that both $v_4$ and $v_5$ are not good, then by $S_2,S_{12}\nsubseteq G$, there exists at least one face $f_k$ ($k\in\{3,5\}$) such that $n_{5^+}(f_k)\geq2$, and thus $ch_1(v)\geq3-3\times\frac{2}{3}-2\times\frac{1}{2}=0$ by R1 and R3. Otherwise if $i\neq j$, then $ch_1(v)\geq3-3\times\frac{2}{3}-2\times\frac{1}{2}=0$ by R1 and R3.

Assume that $d(f_3)=4$. If $d(v_5)=5$, then $ch_1(v)\geq3-3\times\frac{2}{3}-2\times\frac{1}{2}=0$ by R1 and R3. Next we discuss $d(v_5)=4$, and
recall that $v_5\in N(v)$ is not good. On the other hand, since $S_3\nsubseteq G$, there exists at least one face $f_i$ in $\{f_4,f_5\}$ such that $n_{5^+}(f_i)\geq2$. We may also assume that $n_{5^+}(f_2)=1$ (otherwise $ch_1(v)\geq0$). If $n_{5^+}(f_4)\geq2$, then by $S_{11}\nsubseteq G$, $d(v_1)\geq6$, and thus $ch_1(v)\geq3-3\times\frac{2}{3}-\frac{1}{2}-\frac{4}{9}=\frac{1}{18}>0$ by R1 and R3. If $n_{5^+}(f_5)\geq2$, then we are going to claim that $n_{5^+}(f_5)\geq3$. Note that $d(v_3)=d(v_4)=d(v_5)=4$, and we may assume none of them is rich (otherwise $ch_2(v)\geq0$). By $S_2\nsubseteq G$, we get $n_{5^+}(f_5)\geq3$, as claimed. Recall that $n_{5^+}(f_3)\geq2$, we get $d(x)\geq5$. On the other hand, by $S_{16}\nsubseteq G$, we get $d(x)\geq6$. Thus $ch_1(v)\geq3-3\times\frac{2}{3}-\frac{5}{9}-\frac{4}{9}=0$ by R1 and R3.
\end{proof}

\begin{claim}\label{W25}
For each vertex $v\in W_2$ with $f_3(v)=0$, $ch_2(v)\ge0$.
\end{claim}
\begin{proof}Assume that $f_4(v)=0$, then by $S_3\nsubseteq G$, there is at least one $f_i$ $(i\in[5])$ satisfying $f_i\neq (4,4,4,4,5)$, and so $\send{v}{f_i}\leq \frac{1}{2}$. Hence, $ch_1(v)\geq3-3\times\frac{2}{3}-\frac{1}{2}=0$ by R1 and R3. Assume that $f_4(v)=1$. W.l.o.g., let $f_1=(v,v_1,x,v_2)$ be the $4$-face.

\textbf{Case 1.} $n_{5^+}(f_1)=1$. Assume that $f_{6^+}(v)=1$. If $d(f_2)\geq6$, then $n_{5^+}(f_5)\geq2$ by $S_{4}\nsubseteq G$, and $n_{5^+}(f_i)\geq2$ for some $i\in\{3,4\}$ by $S_3\nsubseteq G$. Thus $ch_1(v)\geq3-1-\frac{2}{3}-\frac{1}{3}-2\times\frac{1}{2}=0$ by R3. If $d(f_3)\geq6$, then $n_{5^+}(f_2)\geq2$ and $n_{5^+}(f_5)\geq2$ by $S_{4}\nsubseteq G$. Thus $ch_1(v)\geq3-1-\frac{2}{3}-\frac{1}{3}-2\times\frac{1}{2}=0$ by R1 and R3.

Assume that $f_{6^+}(v)=0$. Since $S_{4}\nsubseteq G$, $n_{5^+}(f_2)\geq2$ and $n_{5^+}(f_5)\geq2$ hold. If $d(v_4)\geq5$, then $ch_1(v)\geq3-1-4\times\frac{1}{2}=0$ by R3. If $d(v_1)=4$, then we may assume that $v_4$ is not good (otherwise $ch_2(v)\geq0$). By $S_2,S_{23}\nsubseteq G$, we get that $n_{5^+}(f_3)\geq2$ and $n_{5^+}(f_4)\geq2$, and thus $ch_1(v)\geq3-1-4\times\frac{1}{2}=0$ by R3.

\textbf{Case 2.}  $n_{5^+}(f_1)\geq2$. Then $ch_1(v)\geq3-\frac{2}{3}-\max\{2\times\frac{2}{3}+\frac{1}{2}+\frac{1}{3}, 2\times\frac{2}{3}+2\times\frac{1}{2}\}=0$ by R3.
\end{proof}
According to all above claims, we know that the minimum counterexample does not exist.
\end{proof}

\section{Proof of Theorem \ref{thm}}\label{proof}
Let $G$ be a counterexample to Theorem \ref{thm} with fewest vertices and edges, that is, there is a list assignment $L$ of $G$ satisfying $|L(v)|\geq 4$ for any $v\in V(G)$   such that $G$ is not $L$-colorable but any proper subgraph of $G$ is $L$-colorable. Firstly, we present the well-known Combinatorial Nullstellensatz initiated by Alon which is essential to produce reducible subgraphs.

\begin{lemma}[\cite{Alon}, Combinatorial Nullstellensatz]\label{Alonlem}
Let $F$ be an arbitrary field, and let $f=f(x_1,\ldots, x_n)$ be a polynomial in $F[x_1,\ldots, x_n]$. Suppose the degree $\deg(f)$ of $f$ is $\sum_{i=1}^nt_i$, where each $t_i$ is a nonnegative integer, and suppose the coefficient of $\prod_{i=1}^nx_i^{t_i}$ in $f$ is nonzero. Then, if $C_1,\ldots, C_n$ are subsets of $F$ with $|C_i|>t_i$, there are $c_1\in C_1, c_2\in C_2, \ldots, c_n\in C_n$ so that
$$f(c_1,\ldots, c_n)\neq0.$$
\end{lemma}

If $G$ has a vertex $v$ of degree at most three, then we can extend an $L$-coloring $\varphi$ of $G\setminus v$ to an $L$-coloring $\phi$ of $G$ by setting $\phi(v)\in L(v)\backslash\{\varphi(u): \; uv\in E(G)\}$, a contradiction. So $\delta(G)\geq 4$. By Lemma \ref{struclemma}, $G$ must contain a subgraph isomorphic to one of the configurations in  $\mathcal{S}$ (see Appendix \ref{appb}). Next, we prove that all these subgraphs do not exist, that is, all configurations $S_1$-$S_{47}$ in $\mathcal{S}$ are reducible, which leads to a contradiction.

\begin{lemma}\label{redulem}
$S_1$-$S_{47}$ in $\mathcal{S}$ are reducible.


\end{lemma}

\begin{proof}
By the minimality of $G$, there is an
$L$-coloring of $G-S_i$ for each $i\in [47]$.
Fix some $i$, say $i_0$, there is an $L$-coloring $\varphi$ of $G-S_{i_0}$. Let $S_{i_0}=\{x_0,x_1,\ldots,x_{n-1}\}$ and  $C_{\varphi}(v)=\{\varphi(u): uv\in E(G)\ and \ u\in V(G-S_{i_0})\}$. Let $C_j=L(x_j)\backslash C_{\varphi}(x_j)$ for $j\in\{0,1,\ldots,n-1\}$. Now we extend $\varphi$ to $G$ and let $\phi$
denote the coloring after all vertices in $S_{i_0}$ are colored. Let $c_0,c_1,\ldots,c_{n-1}$ correspond to the colors of
$x_0,x_1,\ldots,x_{n-1}$ respectively.  If $c_i-c_j\neq0$ for any $x_ix_j\in E(G)$, then $\phi$ is a proper $L$-coloring of $G$. Next let $P=P(x_0,x_1,\ldots,x_{n-1})$ be the following polynomial:
\begin{align}
&P(x_0,x_1,\ldots,x_{n-1})=\displaystyle{\prod^{}_{x_ix_j\in E(G)}}(x_i-x_j).\nonumber
\end{align}
That is, if there are $c_0\in C_0, c_1\in C_1,\ldots,c_{n-1}\in C_{n-1}$ such that $P(x_0,x_1,\ldots,x_{n-1})\neq0$, then we can extend $\varphi$ to an $L$-coloring $\phi$ of $G$ by choosing $x_0=c_0, x_1 = c_1,\ldots,x_{n-1}=c_{n-1}$.

Based on Lemma \ref{Alonlem}, we present an algorithm in Appendix \ref{algori} which effectively calculates reducible configurations. Let us take $S_1$ as an example. Let $S_1=\{x_0,x_1,\ldots,x_4\}$ such that $x_0x_1x_4$ is a triangle and  $x_1x_2x_3x_4$ is a 4-face, where $d(x_i)=4$ for each $i\in\{0,1,2,3\}$ and $d(x_4)=5$. Then $$P(x_0,x_1,\ldots,x_{4})=(x_0-x_1)(x_0-x_4)(x_1-x_2)(x_1-x_4)(x_2-x_3)(x_3-x_4).$$
That is, input ``vve = [(0, 1), (0, 4), (1, 2), (1, 4), (2, 3), (3, 4)]''.
Note that  $|C_1|>2$ and $|C_i|>1$ for each $i\in\{0,2,3,4\}$ as $x_1$ has one neighbor in $V(G-S_1)$ and each $x_i$ has two neighbors in $V(G-S_1)$.  Thus, we
input ``v\_List = [1,2,1,1,1]''. Through  the computation of the algorithm in Appendix \ref{algori},
we get the 1st valid expansion is [1,2,1,1,1], that is,  the coefficient of $x_0x_1^2x_2x_3x_4$ in $P$ is nonzero. Therefore, $S_1$ is reducible by Lemma  \ref{Alonlem}.
\end{proof}

This completes the proof.

\section*{Acknowledgements}
The authors would like to thank the reviewers for their valuable comments, which greatly improve the paper.

\begin{appendix}
\section{Algorithm}\label{algori}
\footnotesize
\begin{lstlisting}
# -*- coding: utf-8 -*-
#!/usr/bin/env python
import copy

def choosable(n,v_List,edges): # Determine whether satisfying \
Combinatorial Nullstellensatz, back to the remainder of the expansion!
    # n: the number of vertices, v_List[0..n-1]: |L(v)|-1, edges: |L(e)|
    zks={}
    zks['0'*n]=1
    len_edges=len(edges)
    for i in range(len_edges):
        v1,v2=edges[i]
        List_zks=[]
        while zks:
            List_zks.append(zks.popitem())
        while List_zks:
            a,b=List_zks.pop()
            if ord(a[v1])-ord("0")<v_List[v1]:
                a1=a[:v1]+chr(ord(a[v1])+1)+a[v1+1:]
                if a1 in zks.keys():
                    zks[a1]=zks[a1]+b
                    if zks[a1]==0:
                        del zks[a1]
                else:
                    zks[a1]=b
            if ord(a[v2])-ord("0")<v_List[v2]:
                a2=a[:v2]+chr(ord(a[v2])+1)+a[v2+1:]
                if a2 in zks.keys():
                    zks[a2]=zks[a2]-b
                    if zks[a2]==0:
                        del zks[a2]
                else:
                    zks[a2]=-b
    return zks

# The main program
def Comb_Null(vve, v_List):
    # ==========================================================================
    # List coloring.
    # vve: Labelling vertices must start at 0. \
    e.g. 3-cycle: vve=[(0, 1),(1, 2),(2, 0)]
    # v_List: |L(v)|-1, must be integers. \
    e.g. 3-cycle: v_List=[1, 1, 1]

    #  Apply Combinatorics Nullstellensatz =======================================================
    v_no=len(v_List)
    zks=choosable(v_no,v_List, vve)

    #  Output part. If there are too many expansions that \
    satisfy the criteria, we print up to 10 =======================
    size_zks=len(zks)
    if size_zks>0:
        print("\n\nThe total number of valid expansions= "+str(size_zks)+",\
        among them:")
        i=percent = 0
        for a in zks.keys():
            if i/size_zks>=percent:
                if i == 0:
                    print("The 1st valid expansion is: [",end="")
                elif i == 1:
                    print("The 2nd valid expansion is: [",end="")
                else:
                    print("The "+str(i+1)+"th valid expansion is: [",end="")
                for j in range(v_no-1):
                    print(str(ord(a[j])-ord("0"))+",",end="")
                print(str(ord(a[v_no-1])-ord("0"))+"]")
                percent+=0.1
            i+=1
    else: print("\n\n No valid expansion!!")

#Example
#Input
vve = [(0, 1), (0, 4), (1, 2), (1, 4), (2, 3), (3, 4)]
v_List = [1, 2, 1, 1, 1]
Comb_Null(vve, v_List)

# Output
# The total number of valid expansions= 1, among them:
# The 1st valid expansion is: [1,2,1,1,1]
\end{lstlisting}

\section{All configurations in $\mathcal{S}$}\label{appb}
\begin{figure}[H]
 \begin{center}
   \includegraphics[scale=0.67]{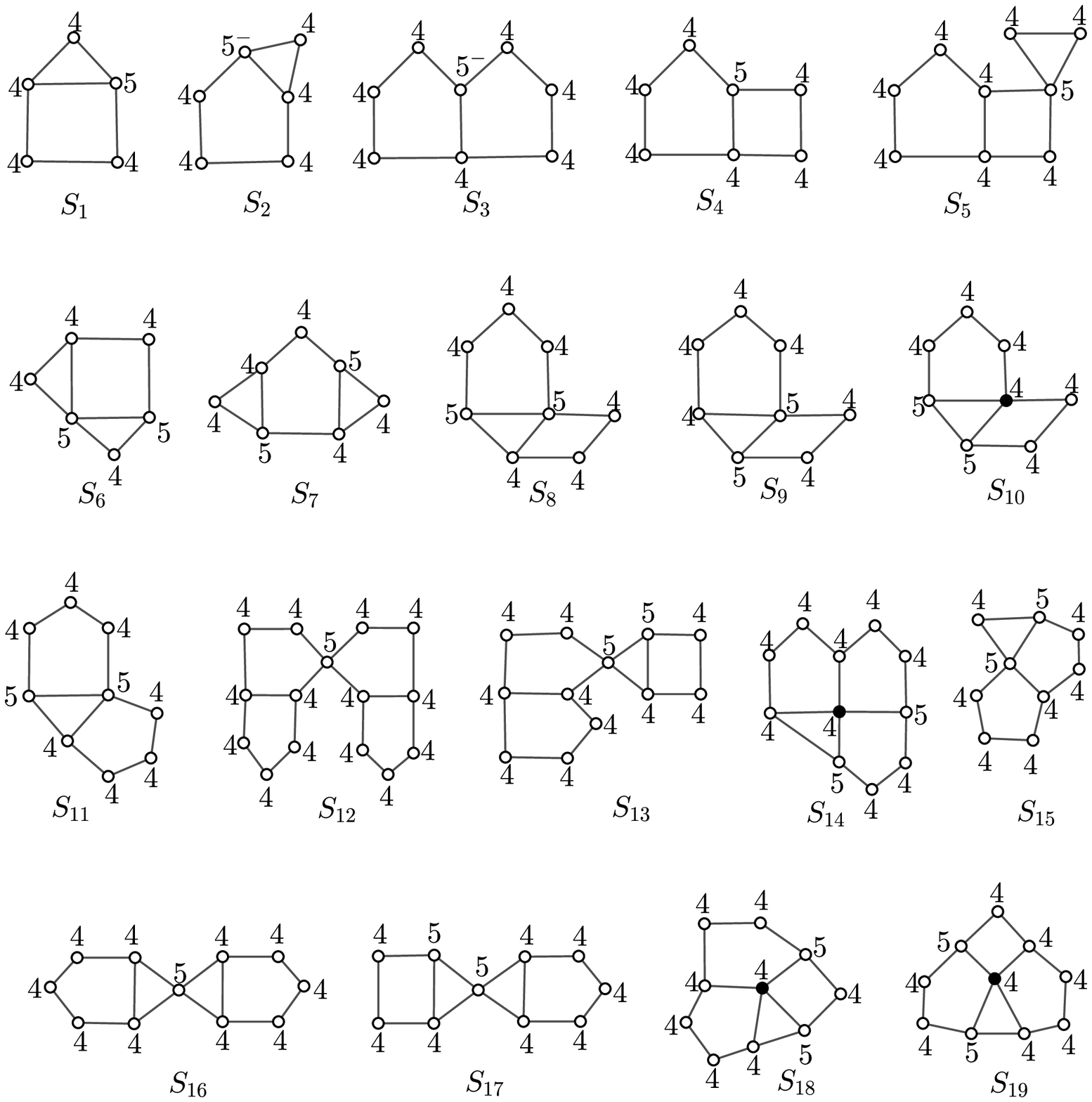}\\
   \label{redu1fig}
 \end{center}
\end{figure}
\begin{figure}[H]
 \begin{center}
   \includegraphics[scale=0.67]{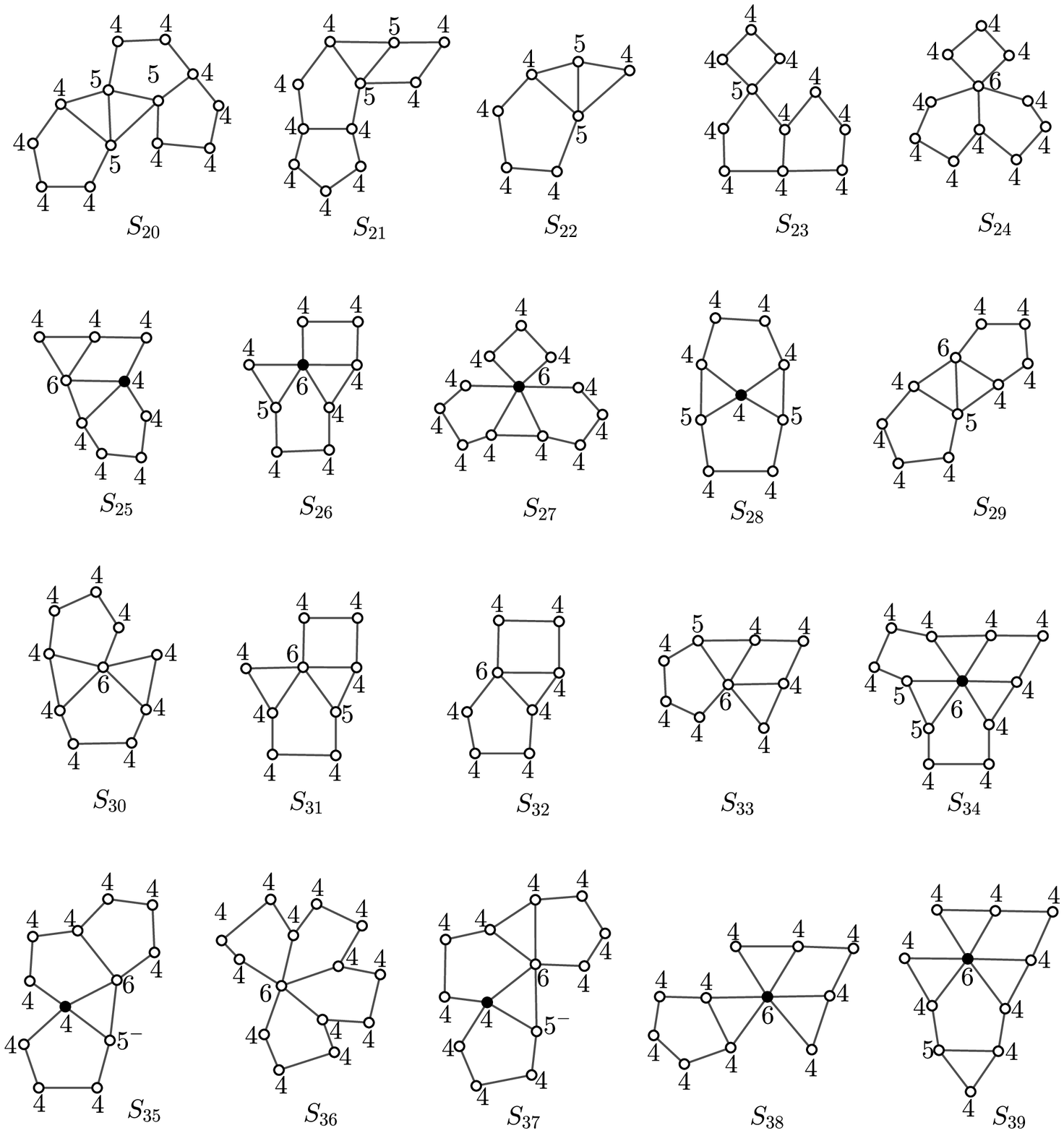}\\
   \label{redu1fig}
 \end{center}
\end{figure}
\vspace{-1cm}
\begin{figure}[H]
 \begin{center}
   \includegraphics[scale=0.67]{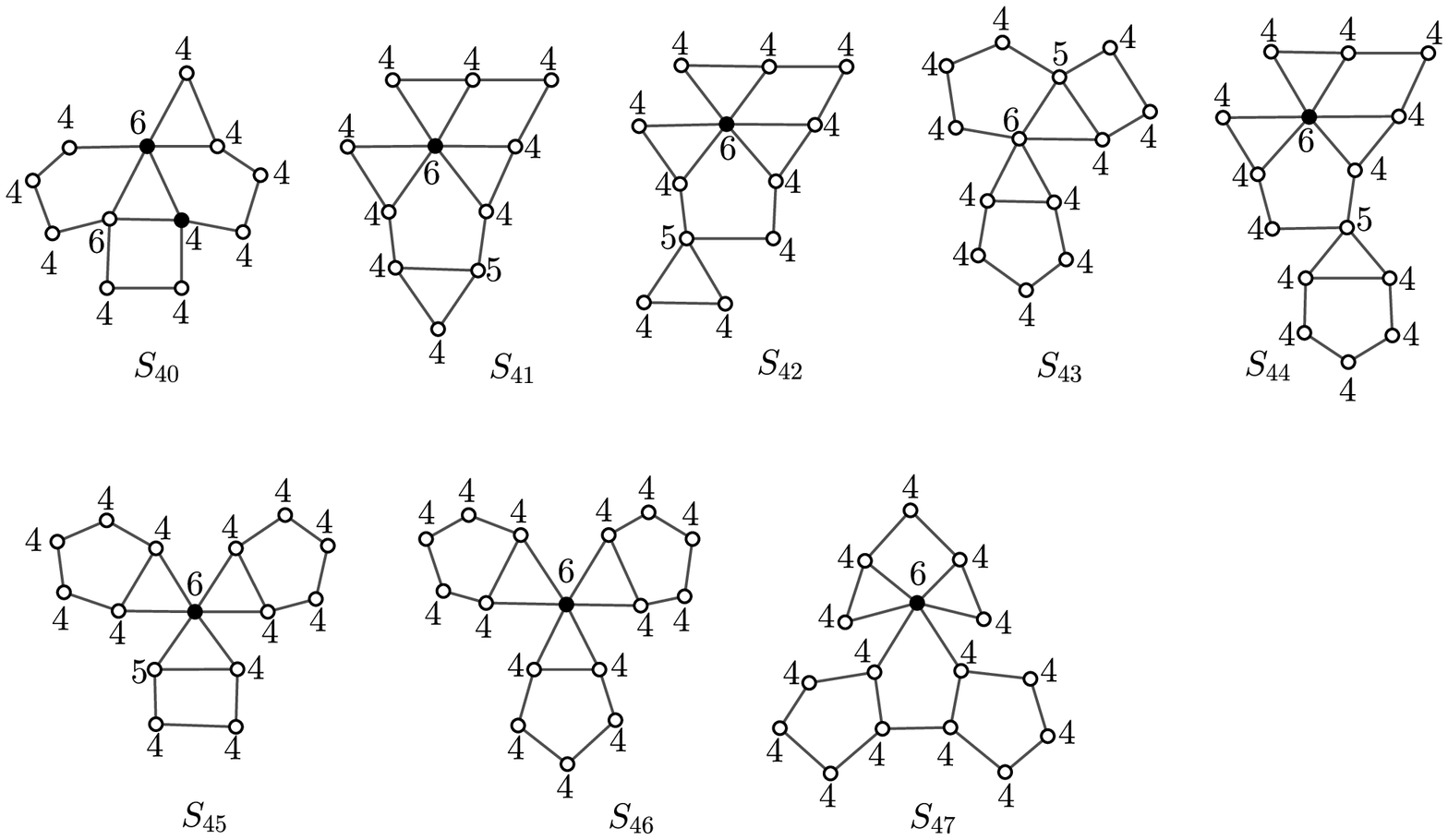}
 \end{center}
\end{figure}

\end{appendix}
\end{document}